\magnification \magstep1


\ifx\begin\undefined\else\global\advance\srcdepth by
1\expandafter \fi

\def\begin{}
\newcount\srcdepth
\srcdepth=1

\outer\def\bye{\global\advance\srcdepth by -1
  \ifnum\srcdepth=0
    \def\endcmd{\vfill\eject\nopagenumbers\par\vfill\supereject\end}
  \else\def\endcmd{}\fi
  \endcmd
}


\baselineskip=13pt
\hsize = 5.5truein
\hoffset = 0.5truein
\vsize = 8.5truein
\voffset = 0.2truein
\emergencystretch = 0.05\hsize

\overfullrule=0pt

\newif\ifblackboardbold

\blackboardboldtrue


\font\sectionfont=cmbx12


\newfam\bboldfam
\ifblackboardbold
\font\tenbbold=msbm10
\font\sevenbbold=msbm7
\font\fivebbold=msbm5
\textfont\bboldfam=\tenbbold
\scriptfont\bboldfam=\sevenbbold
\scriptscriptfont\bboldfam=\fivebbold
\def\bbold{\fam\bboldfam\tenbbold}
\else
\def\bbold{\bf}
\fi


\font\Arm=cmr8
\font\Ai=cmmi8
\font\Asy=cmsy8
\font\Abf=cmbx8
\font\Brm=cmr6
\font\Bi=cmmi6
\font\Bsy=cmsy6
\font\Bbf=cmbx6
\font\Crm=cmr5
\font\Ci=cmmi5
\font\Csy=cmsy5
\font\Cbf=cmbx5

\ifblackboardbold
\font\Abbold=msbm10 at 8pt
\font\Bbbold=msbm7 at 6pt
\font\Cbbold=msbm5
\fi

\def\smallmath{%
\textfont0=\Arm \scriptfont0=\Brm \scriptscriptfont0=\Crm
\textfont1=\Ai \scriptfont1=\Bi \scriptscriptfont1=\Ci
\textfont2=\Asy \scriptfont2=\Bsy \scriptscriptfont2=\Csy
\textfont\bffam=\Abf \scriptfont\bffam=\Bbf \scriptscriptfont\bffam=\Cbf
\def\rm{\fam0\Arm}\def\mit{\fam1}\def\oldstyle{\fam1\Ai}%
\def\bf{\fam\bffam\Abf}%
\ifblackboardbold
\textfont\bboldfam=\Abbold
\scriptfont\bboldfam=\Bbbold
\scriptscriptfont\bboldfam=\Cbbold
\def\bbold{\fam\bboldfam\Abbold}%
\fi
}








\newlinechar=`@
\def\forwardmsg#1#2#3{\immediate\write16{@*!*!*!* forward reference should
be: @\noexpand\forward{#1}{#2}{#3}@}}
\def\nodefmsg#1{\immediate\write16{@*!*!*!* #1 is an undefined reference@}}

\def\forwardsub#1#2{\def\newref{{#2}{#1}}}

\def\forward#1#2#3{%
\expandafter\expandafter\expandafter\forwardsub\expandafter{#3}{#2}
\expandafter\ifx\csname#1\endcsname\relax\else%
\expandafter\ifx\csname#1\endcsname\newref\else%
\forwardmsg{#1}{#2}{#3}\fi\fi%
\expandafter\let\csname#1\endcsname\newref}

\def\firstarg#1{\expandafter\argone #1}\def\argone#1#2{#1}
\def\secondarg#1{\expandafter\argtwo #1}\def\argtwo#1#2{#2}

\def\ref#1{\expandafter\ifx\csname#1\endcsname\relax
  {\nodefmsg{#1}\bf`#1'}\else
  \expandafter\firstarg\csname#1\endcsname
  ~\expandafter\secondarg\csname#1\endcsname\fi}

\def\refs#1{\expandafter\ifx\csname#1\endcsname\relax
  {\nodefmsg{#1}\bf`#1'}\else
  \expandafter\firstarg\csname #1\endcsname
  s~\expandafter\secondarg\csname#1\endcsname\fi}

\def\refn#1{\expandafter\ifx\csname#1\endcsname\relax
  {\nodefmsg{#1}\bf`#1'}\else
  \expandafter\secondarg\csname #1\endcsname\fi}



\def\widow#1{\vskip 0pt plus#1\vsize\goodbreak\vskip 0pt plus-#1\vsize}



\def\marginlabel#1{}

\def\showlabelsabove{
\font\labelfont=cmss10 at 6pt
\def\marginlabel##1{\rlap{\smash{\raise 10pt\hbox{\labelfont##1}}}}
}

\newcount\seccount
\newcount\proccount
\seccount=0
\proccount=0

\def\stdskip{\vskip 9pt plus3pt minus 3pt}
\def\stdbreak{\par\removelastskip\penalty-100\stdskip}

\def\proof{\stdbreak\noindent{\sl Proof. }}

\def\qed{\vrule height 1.2ex width .9ex depth .1ex}

\def\Box{
  \ifmmode\eqno\qed
  \else\ifvmode\removelastskip\line{\hfil\qed}
  \else\unskip\quad\hskip-\hsize
    \hbox{}\hskip\hsize minus 1em\qed\par
  \fi\stdbreak\fi}

\def\references{
  \removelastskip
  \widow{.05}
  \vskip 24pt plus 6pt minus 6 pt
  \leftline{\sectionfont References}
  \nobreak\stdskip\noindent}

\def\ifempty#1#2\endB{\ifx#1\endA}
\def\makeref#1#2#3{\ifempty#1\endA\endB\else\forward{#1}{#2}{#3}\fi}

\outer\def\section#1 #2\par{
  \removelastskip
  \global\advance\seccount by 1
  \global\proccount=0\relax
                \edef\numtoks{\number\seccount}
  \makeref{#1}{Section}{\numtoks}
  \widow{.05}
  \vskip 24pt plus 6pt minus 6 pt
  \message{#2}
  \leftline{\marginlabel{#1}\sectionfont\numtoks\quad #2}
  \nobreak\stdskip}

\def\proclamation#1#2{
  \outer\expandafter\def\csname#1\endcsname##1 ##2\par{
  \stdbreak
  \advance\proccount by 1
  \edef\numtoks{\number\seccount.\number\proccount}
  \makeref{##1}{#2}{\numtoks}
  \noindent{\marginlabel{##1}\bf #2 \numtoks\enspace}
  {\sl##2\par}
  \stdbreak}}

\def\othernumbered#1#2{
  \outer\expandafter\def\csname#1\endcsname##1{
  \stdbreak
  \advance\proccount by 1
  \edef\numtoks{\number\seccount.\number\proccount}
  \makeref{##1}{#2}{\numtoks}
  \noindent{\marginlabel{##1}\bf #2 \numtoks\enspace}}}

\proclamation{definition}{Definition}
\proclamation{lemma}{Lemma}
\proclamation{proposition}{Proposition}
\proclamation{theorem}{Theorem}
\proclamation{corollary}{Corollary}
\proclamation{conjecture}{Conjecture}

\othernumbered{example}{Example}
\othernumbered{remark}{Remark}
\othernumbered{construction}{Construction}






\newcount\figcount
\figcount=0
\newbox\drawing
\newcount\drawbp
\newdimen\drawx
\newdimen\drawy
\newdimen\ngap
\newdimen\sgap
\newdimen\wgap
\newdimen\egap

\def\drawbox#1#2#3{\vbox{
  \setbox\drawing=\vbox{\offinterlineskip\epsfbox{#2.eps}\kern 0pt}
  \drawbp=\epsfurx
  \advance\drawbp by-\epsfllx\relax
  \multiply\drawbp by #1
  \divide\drawbp by 100
  \drawx=\drawbp truebp
  \ifdim\drawx>\hsize\drawx=\hsize\fi
  \epsfxsize=\drawx
  \setbox\drawing=\vbox{\offinterlineskip\epsfbox{#2.eps}\kern 0pt}
  \drawx=\wd\drawing
  \drawy=\ht\drawing
  \ngap=0pt \sgap=0pt \wgap=0pt \egap=0pt 
  \setbox0=\vbox{\offinterlineskip
    \box\drawing \ifgridlines\drawgrid\drawx\drawy\fi #3}
  \kern\ngap\hbox{\kern\wgap\box0\kern\egap}\kern\sgap}}

\def\draw#1#2#3{
  \setbox\drawing=\drawbox{#1}{#2}{#3}
  \advance\figcount by 1
  \goodbreak
  \midinsert
  \centerline{\ifgridlines\boxgrid\drawing\fi\box\drawing}
  \smallskip
  \vbox{\offinterlineskip
    \centerline{Figure~\number\figcount}
    \smash{\marginlabel{#2}}}
  \endinsert}

\def\nextfigtoks{%
  \advance\figcount by 1%
  \edef\numtoks{\number\figcount}%
  \advance\figcount by -1}

\newif\ifgridlines
\newbox\figtbox
\newbox\figgbox
\newdimen\figtx
\newdimen\figty

\newdimen\bwd
\bwd=2sp 

\def\hline#1{\vbox{\smash{\hbox to #1{\leaders\hrule height \bwd\hfil}}}}

\def\vline#1{\hbox to 0pt{%
  \hss\vbox to #1{\leaders\vrule width \bwd\vfil}\hss}}

\def\clap#1{\hbox to 0pt{\hss#1\hss}}
\def\vclap#1{\vbox to 0pt{\offinterlineskip\vss#1\vss}}

\def\hstutter#1#2{\hbox{%
  \setbox0=\hbox{#1}%
  \hbox to #2\wd0{\leaders\box0\hfil}}}

\def\vstutter#1#2{\vbox{
  \setbox0=\vbox{\offinterlineskip #1}
  \dp0=0pt
  \vbox to #2\ht0{\leaders\box0\vfil}}}

\def\crosshairs#1#2{
  \dimen1=.002\drawx
  \dimen2=.002\drawy
  \ifdim\dimen1<\dimen2\dimen3\dimen1\else\dimen3\dimen2\fi
  \setbox1=\vclap{\vline{2\dimen3}}
  \setbox2=\clap{\hline{2\dimen3}}
  \setbox3=\hstutter{\kern\dimen1\box1}{4}
  \setbox4=\vstutter{\kern\dimen2\box2}{4}
  \setbox1=\vclap{\vline{4\dimen3}}
  \setbox2=\clap{\hline{4\dimen3}}
  \setbox5=\clap{\copy1\hstutter{\box3\kern\dimen1\box1}{6}}
  \setbox6=\vclap{\copy2\vstutter{\box4\kern\dimen2\box2}{6}}
  \setbox1=\vbox{\offinterlineskip\box5\box6}
  \smash{\vbox to #2{\hbox to #1{\hss\box1}\vss}}}

\def\boxgrid#1{\rlap{\vbox{\offinterlineskip
  \setbox0=\hline{\wd#1}
  \setbox1=\vline{\ht#1}
  \smash{\vbox to \ht#1{\offinterlineskip\copy0\vfil\box0}}
  \smash{\vbox{\hbox to \wd#1{\copy1\hfil\box1}}}}}}

\def\drawgrid#1#2{\vbox{\offinterlineskip
  \dimen0=\drawx
  \dimen1=\drawy
  \divide\dimen0 by 10
  \divide\dimen1 by 10
  \setbox0=\hline\drawx
  \setbox1=\vline\drawy
  \smash{\vbox{\offinterlineskip
    \copy0\vstutter{\kern\dimen1\box0}{10}}}
  \smash{\hbox{\copy1\hstutter{\kern\dimen0\box1}{10}}}}}

\def\figtext#1#2#3#4#5{
  \setbox\figtbox=\hbox{#5}
  \dp\figtbox=0pt
  \figtx=-#3\wd\figtbox \figty=-#4\ht\figtbox
  \advance\figtx by #1\drawx \advance\figty by #2\drawy
  \dimen0=\figtx \advance\dimen0 by\wd\figtbox \advance\dimen0 by-\drawx
  \ifdim\dimen0>\egap\global\egap=\dimen0\fi
  \dimen0=\figty \advance\dimen0 by\ht\figtbox \advance\dimen0 by-\drawy
  \ifdim\dimen0>\ngap\global\ngap=\dimen0\fi
  \dimen0=-\figtx
  \ifdim\dimen0>\wgap\global\wgap=\dimen0\fi
  \dimen0=-\figty
  \ifdim\dimen0>\sgap\global\sgap=\dimen0\fi
  \smash{\rlap{\vbox{\offinterlineskip
    \hbox{\hbox to \figtx{}\ifgridlines\boxgrid\figtbox\fi\box\figtbox}
    \vbox to \figty{}
    \ifgridlines\crosshairs{#1\drawx}{#2\drawy}\fi
    \kern 0pt}}}}


\def\hpad#1#2#3{\hbox{\kern #1\hbox{#3}\kern #2}}
\def\vpad#1#2#3{\setbox0=\hbox{#3}\dp0=0pt\vbox{\kern #1\box0\kern #2}}



\def\stack#1#2#3{\vbox{\offinterlineskip
  \setbox2=\hbox{#2}
  \setbox3=\hbox{#3}
  \dimen0=\ifdim\wd2>\wd3\wd2\else\wd3\fi
  \hbox to \dimen0{\hss\box2\hss}
  \kern #1
  \hbox to \dimen0{\hss\box3\hss}}}


\def\hexp#1{%
  \setbox0=\hbox{${}^{#1}$}%
  \hbox to .5\wd0{\box0\hss}}



\def\bmatrix#1#2{{\smallmath\left[\vcenter{\halign
  {&\kern#1\hfil$##\mathstrut$\kern#1\cr#2}}\right]}}

\def\rightarrowmat#1#2#3{
  \setbox1=\hbox{\kern#2$\bmatrix{#1}{#3}$\kern#2}
  \,\vbox{\offinterlineskip\hbox to\wd1{\hfil\copy1\hfil}
    \kern 3pt\hbox to\wd1{\rightarrowfill}}\,}

\def\leftarrowmat#1#2#3{
  \setbox1=\hbox{\kern#2$\bmatrix{#1}{#3}$\kern#2}
  \,\vbox{\offinterlineskip\hbox to\wd1{\hfil\copy1\hfil}
    \kern 3pt\hbox to\wd1{\leftarrowfill}}\,}

\def\rightarrowbox#1#2{
  \setbox1=\hbox{\kern#1\hbox{\smallmath #2}\kern#1}
  \,\vbox{\offinterlineskip\hbox to\wd1{\hfil\copy1\hfil}
    \kern 3pt\hbox to\wd1{\rightarrowfill}}\,}

\def\leftarrowbox#1#2{
  \setbox1=\hbox{\kern#1\hbox{\smallmath #2}\kern#1}
  \,\vbox{\offinterlineskip\hbox to\wd1{\hfil\copy1\hfil}
    \kern 3pt\hbox to\wd1{\leftarrowfill}}\,}






\def\bookletdims{
  \hsize=5.25truein
  \vsize=7truein
}

\def\legalbooklet#1{
  \input quire
  \bookletdims
  \htotal=7.0truein
  \vtotal=8.5truein
  \hoffset=\htotal
  \advance\hoffset by -\hsize
  \divide\hoffset by 2
  \voffset=\vtotal
  \advance\voffset by -\vsize
  \divide\voffset by 2
  \advance\voffset by -.0625truein
  \shhtotal=2\htotal
  \horigin=0.0truein
  \vorigin=0.0truein
  \shstaplewidth=0.01pt
  \shstaplelength=0.66truein
  \shthickness=0pt
  \shoutline=0pt
  \shcrop=0pt
  \shvoffset=-1.0truein
  \ifnum#1>0\quire{#1}\else\qtwopages\fi
}

\def\preview{
  \input quire
  \bookletdims
  \hoffset=0.1truein
  \vtotal=8.5truein
  \shhtotal=14truein
  \voffset=\vtotal
  \advance\voffset by -\vsize
  \divide\voffset by 2
  \advance\voffset by -.0625truein
  \htotal=2\hoffset
  \advance\htotal by \hsize
  \horigin=0.0truein
  \vorigin=0.0truein
  \shstaplewidth=0.5pt
  \shstaplelength=0.5\vtotal
  \shthickness=0pt
  \shoutline=0pt
  \shcrop=0pt
  \shvoffset=-1.0truein
  \qtwopages
}

\def\twoup{
  \input quire
  \hsize=4.79452truein 
  \vsize=7truein
  \vtotal=8.5truein
  \shhtotal=11truein
  \hoffset=-2\hsize
  \advance\hoffset by \shhtotal
  \divide\hoffset by 6
  \voffset=\vtotal
  \advance\voffset by -\vsize
  \divide\voffset by 2
  \advance\voffset by -12truept
  \htotal=2\hoffset
  \advance\htotal by \hsize
  \horigin=0.0truein
  \vorigin=0.0truein
  \shstaplewidth=0.01pt
  \shstaplelength=0pt
  \shthickness=0pt
  \shoutline=0pt
  \shcrop=0pt
  \shvoffset=-1.0truein
  \qtwopages
}


\newcount\countA
\newcount\countB
\newcount\countC

\def\monthname{\begingroup
  \ifcase\number\month
    \or January\or February\or March\or April\or May\or June\or
    July\or August\or September\or October\or November\or December\fi
\endgroup}

\def\dayname{\begingroup
  \countA=\number\day
  \countB=\number\year
  \advance\countA by 0 
  \advance\countA by \ifcase\month\or
    0\or 31\or 59\or 90\or 120\or 151\or
    181\or 212\or 243\or 273\or 304\or 334\fi
  \advance\countB by -1995
  \multiply\countB by 365
  \advance\countA by \countB
  \countB=\countA
  \divide\countB by 7
  \multiply\countB by 7
  \advance\countA by -\countB
  \advance\countA by 1
  \ifcase\countA\or Sunday\or Monday\or Tuesday\or Wednesday\or
    Thursday\or Friday\or Saturday\fi
\endgroup}

\def\timename{\begingroup
   \countA = \time
   \divide\countA by 60
   \countB = \countA
   \countC = \time
   \multiply\countA by 60
   \advance\countC by -\countA
   \ifnum\countC<10\toks1={0}\else\toks1={}\fi
   \ifnum\countB<12 \toks0={\sevenrm AM}
     \else\toks0={\sevenrm PM}\advance\countB by -12\fi
   \relax\ifnum\countB=0\countB=12\fi
   \hbox{\the\countB:\the\toks1 \the\countC \thinspace \the\toks0}
\endgroup}

\def\timestamp{\dayname, \the\day\ \monthname\ \the\year, \timename}


\def\enma#1{{\ifmmode#1\else$#1$\fi}}

\input diagrams.tex
\def\A{{\cal A}}

\def \I{{\cal I}}
\def\R{{\cal R}}
\def\O{{\cal O}}

\def\FF{{\bf F}}
\def\GG{{\bf G}}
\def\KK{{\bf K}}
\def\LL{{\bf L}}
\def\II{{\bf I}}

\def\mm{{\bf m}}

\def\mindeg{\mathop{\rm mindeg}\nolimits}
\def\depth{\mathop{\rm depth}\nolimits}
\def\In{\mathop{\rm in}\nolimits}
\def\Sym{\mathop{\rm Sym}\nolimits}

\def\H{\mathop{\rm H}\nolimits}
\def\Tor{\mathop{\rm Tor}\nolimits}
\def\Ext{\mathop{\rm Ext}\nolimits}
\def\Hom{\mathop{\rm Hom}\nolimits}

\def\PP{\mathop{\bf P}\nolimits}
\def\rank{\mathop{\rm rank}\nolimits}
\def\reg{\mathop{\rm reg}\nolimits}
\def\codim{\mathop{\rm codim}\nolimits}

\def\Tot{\mathop{\rm Tot}\nolimits}
\def\In{\mathop{\rm in}\nolimits}
\def\gin{\mathop{\rm Gin}\nolimits}
\def\coker{\mathop{\rm coker}\nolimits}

\forward {deg syz}{Section}{2}
\forward {C-M reg section}{Section}{3}
\forward {subadd section}{Section}{4}
\forward {specializing section}{Section}{5}
\forward {quadrics section}{Section}{6}
\forward {reg prods section}{Section}{7}
\forward {monomial ideals}{Section}{8}
\forward {symmetric algebra section}{Section}{9}
\forward {elimination section}{Section}{10}
\forward {almost lin}{Section}{11}

\forward {truncation}{Proposition}{1.7}
\forward {nonstandard grading}{Remark}{2.4}
\forward{sliding}{Lemma}{2.6}
\forward {sharpness}{Proposition}{2.7}
\forward{reg of Tor}{Corollary}{3.1}
\forward{regularity of general sections}{Corollary}{3.9} 
\forward{saturation of plane sections}{Corollary}{3.13}
\forward {socle estimation}{Corollary}{4.2}
\forward{half-way linear2}{Corollary}{7.6}
\forward{reg of powers}{Corollary}{7.8}
\forward{reg of powers2}{Corollary}{7.9}
\forward {non annihilation}{Example}{9.3}

\let\iso\cong
 
\def\emph#1{{\it#1\/}} 
\centerline{\bf The Regularity of Tor and Graded Betti Numbers}
\smallskip
\centerline{by}
\smallskip
\centerline{{David Eisenbud, Craig Huneke and Bernd Ulrich}
\footnote{$^\dagger$}{\noindent We wish to thank MSRI, 
where we were all guests (even Eisenbud)
during most of the work on this paper; and to the NSF for its
generous support.}
\footnote{}{Revised 5/18/04}}
\bigskip

\noindent{\bf Abstract:} { 
Let $S=K[x_1,\dots, x_n]$, let $A,B$ be finitely
generated graded $S$-modules, and let $\mm=(x_1,\dots,x_n)\subset S$.
We give bounds for the regularity
of the local cohomology of $\Tor_k(A,B)$
in terms of the graded Betti numbers of $A$ and $B$,
under the assumption that $\dim\Tor_1(A,B) \leq 1.$
We apply the results to syzygies,
Gr\"obner bases, 
products and powers of ideals,
and to the relationship of the Rees
and Symmetric algebras. For example we show that
any homogeneous linearly presented $\mm$-primary
ideal has some power equal to a power of $\mm$; and if 
the first $\lceil(n-1)/2\rceil$ steps
 of the resolution of $I$ are linear, then
$I^2$ is a power of $\mm$. 
 }

\section{intro} Introduction

Let $S=K[x_1,\dots, x_n]$ and let $A,B$ be finitely
generated graded $S$-modules. If $T$ is a finitely generated
graded $S$-module we write $\reg T$ for the 
Castelnuovo-Mumford regularity of $T$, and we extend this
to Artinian modules $T$ by setting
$\reg T=\max\{i\mid T_i\neq 0\}.$ 
The main technical results of this paper,
proved in \ref{deg syz}, give upper bounds on the regularity
of the local cohomology modules,
$\H^j_\mm(\Tor_k(A,B))$
under the hypothesis that $\Tor_1(A,B)$ has Krull dimension $\leq 1$.
A special case says that if $A\otimes B$ has finite length then,
for any $k$,
$$
\reg \Tor_{k}(A,B) +n
\leq 
\reg \Tor_{p}(A,K) + \reg \Tor_q(B,K)
$$
for any $p,q$ with
$$
p+q=n+k,\ p\leq n,\ q\leq n.
$$ 

In this formula $\Tor_{p}(A,K)$ is just the graded vector
space of generators of the minimal $p$-th syzygies of $A$,
and $\reg \Tor_p(A,K)$ is simply the maximal degree of such
a generator. Such terms occur so often in this paper that we
will adopt a special notation, and write
$$
t_p(A) := \reg \Tor_p(A,K).
$$ 

The rest of the paper is devoted to applications of the
bounds proven in section 2. By way of introduction,
we will now sample the less technical consequences.
Almost every result stated below occurs with more 
generality in the body of the paper.

We begin, in \ref{C-M reg section}, 
with the  regularity of the Tor modules.
We show that if $A$ and $B$ are finitely generated
 graded $S$-modules
such that $\dim \Tor_1(A,B) \leq 1$, then
$$
\reg \Tor_k(A,B) \leq \reg A +\reg B +k,
$$
which generalizes results of Chandler, Sidman, Caviglia
and others. 
For a geometric consequence, let $X, Y\subset \PP^{n-1}$
be projective schemes.
It is elementary that, if $I$ and $J$ are their homogeneous ideals, then
the ideal of forms vanishing on $X\cap Y$ is
equal to $I+J$ in degree $d\gg 0$. It follows from our results
that if $\dim(X\cap Y)=0$ then
it suffices to take
$$
d> t_p(S/I) + t_q(S/J) -n 
$$
for any $p,q$ such that
$p\leq \codim X,\ q\leq \codim Y,$ and $p+q=n.$

In \ref{subadd section} we deduce relations between graded Betti numbers.
For example, we show that if $A=B=S/I$
is a cyclic module of dimension $\leq 1$,
then the function $p\mapsto t_p(S/I)$ satisfies the 
weak convexity condition
$$
 t_n(S/I) \leq  t_{p}(S/I) +  t_{n-p}(S/I)
$$
for $0\leq p\leq n$.

We also compare the graded Betti numbers of a module and an
ideal that annihilates it. We prove that if $S/I$ is
Cohen-Macaulay of codimension $c$, and $I$ contains
a regular sequence of elements of degrees $d_1\leq \cdots \leq d_q$,
then
$$
 t_{c}(S/I) \leq  t_{c-q}(S/I)+d_1+\cdots +d_q.
$$
If $I$ is generated in degrees $\leq d$, then we can take all the
$d_i=d$, and we see that 
$$
 t_c(S/I) -  t_{c-q}(S/I) \leq qd.
$$

In \ref{specializing section} we study the
relationship between the graded Betti numbers 
of an ideal $I$
and its initial ideal in reverse lexicographic order. 
For example, suppose
that
$I\subset S$ is a homogeneous $\mm$-primary ideal generated in 
degree $d$. Setting $m= t_p(S/I)$,
we show that
the initial ideal of $I$ in reverse lexicographic
order contains $(x_1,\dots,x_{p})^{m-p+1}$.
If
the minimal free resolution of $I$ is linear for
$q$ steps and $L$ is any ideal generated by $n-q-1$
independent linear forms, then we show that
$$
\mm^d\subset I+L.
$$
In terms of free resolutions this says $\reg(I+L)\leq d$.

In \ref{quadrics section} we explore the meaning of this
last condition
by characterizing
the ideals $I$ generated by quadrics such that
$
\mm^2\subset I+L
$
for every ideal $L$ generated by $n-q-1$ independent
linear forms.

In \ref{reg prods section} we study powers of linearly presented ideals.
The following conjecture sparked this entire paper:

\conjecture{Eisenbud-Ulrich} (Eisenbud and Ulrich)
If $I\subset S$ is a linearly presented $\mm$-primary
ideal generated in degree $d$, then $I^{n-1}=\mm^{(n-1)d}.$

We prove this conjecture when $n=3$,
and, in \ref{monomial ideals}, for the case of monomial ideals.
But in general we can prove only an asymptotic statement:

\theorem{some power-intro} If $I$ is an $\mm$-primary linearly
presented ideal generated in degree $d$,
then $I^t=\mm^{dt}$ for all $t\gg 0$.
 
The case $n=3$ is generalized by the following result, which
is perhaps the most surprising result of this paper:

\theorem{half-way linear2-intro} Suppose $I$ and $J$ 
are homogeneous ideals
in $S$ of dimension $\leq 1$, generated in degree $d$.
If the resolutions of $I$ and $J$ are linear for 
$\lceil (n-1)/2\rceil$ steps
$($for instance if $I$ and $J$ have linear presentation 
and $n\leq 3$$)$, then
$IJ$ has linear resolution.
In particular, if $I$ and $J$ are $\mm$-primary then $IJ = \mm^{2d}$.

Here the last statement follows from the previous
one because the powers of the maximal ideal are the only $\mm$-primary
ideals with linear resolutions.
Based on this result we were led to generalize \ref{Eisenbud-Ulrich}
as follows:

\conjecture{Us-reg} Suppose that $I$ is an $\mm$-primary
ideal.
\item{$(a)$} If $m\leq kj$, 
then 
$$
t_{m}(I^k)\leq k\; t_j(I) - (kj-m).
$$
\item{$(b)$} If $s\leq n-1$, then 
$$
t_{n-1}(I^k) \leq t_s(I^k)+(n-1-s)t_0(I).
$$
\noindent In particular, taking $s=kj\leq n-1$,
$$
\reg I^k\leq k\; t_j(I)+(n-1-kj)t_0(I)-(n-1).
$$

These formulas are sharp
for complete intersections of forms of degree $d$. Part
$(b)$ is evident for $s=0$: just replace $I$ by a 
complete intersection contained in $I$.

For example, suppose $I$ is generated in degree $d$
and has linear resolution for $j$ 
steps, so $t_j(I)=d+j$.
Taking $k=\lceil(n-1)/j\rceil$, \ref{Us-reg} would give
$\reg I^k\leq kd$, so that $\mm^{kd}\subset I^k$.

Part $(a)$ would imply the truth of the following, which
is in fact equivalent (see \ref{truncation}):
\conjecture{Us} If $I$ is an $\mm$-primary ideal,
and $I$ has a linear resolution for $s$ steps,
then $I^t$ has linear resolution for $st$ steps, 
and $I^t$ is equal to a power of $\mm$
for all $t\geq (n-1)/s$.

We don't even know that powers of $\mm$-primary, linearly
presented ideals are linearly presented!
Examples of Sturmfels [2000] (see also Conca [2003]) show that this
would not be the case without the $\mm$-primary
hypothesis.

\medskip

The torsion in $I\otimes I^t$ is $\Tor_2(S/I, S/I^t)$.
In \ref{symmetric algebra section} we use this relationship
to study the torsion in the symmetric algebra
$\Sym(I)$. We were motivated by the following conjecture of
Eisenbud and Ulrich for linearly presented ideals $I\subset S$
that are of linear type on the punctured spectrum
(that is, each  torsion element in $\Sym I$ is annihilated by
a power of $\mm$): 

\conjecture{sym torsion conj}
If $I$ is generated in degree $d$,
then the torsion in $\Sym_t I$ is generated in
degree $td$.
If  $I$ has linear free resolution,
 then 
the torsion is annihilated by $\mm$; equivalently,
the symmetric algebra of $I$ is a subalgebra of the 
symmetric algebra of the maximal ideal.  

We are able to show, for example, that
if $I$ is an $\mm$-primary ideal generated in degree $d$, 
and has a free resolution that is linear for
$\lceil n/2\rceil$ steps, then, for every $t$, 
the torsion in $\Sym_tI$ is concentrated in degree $dt$.
(Related ideas show that $\wedge^t I$ is
a vector space concentrated in degree $dt$.)
We show in \ref{non annihilation}
 that, at least for $n=3$, the bound $\lceil n/2\rceil$
is sharp. 

In \ref{elimination section} we explore a consequence for elimination 
theory, a method of finding the defining ideal of the 
image of a map $\alpha_V:\PP^{n-1}\to \PP^{N-1}$ 
defined
by an $N$-dimensional vector space $V\subset S_d$ 
of forms of degree $d$.
We assume that the morphism $\alpha_V$ is everywhere defined,
which means that $V$ generates an ideal $I=SV$ that 
is $\mm$-primary. 
Let $M=\dim\Tor_1(I,K)$ 
be the number of
relations required for $I$, and let $\phi$ be the
$N\times M$ matrix of linear forms that presents $I$.
The matrix $\phi$ can be represented as an $n\times N\times M$
tensor over $K$, and thus also represents an $n\times M$ matrix
of linear forms $\psi$ over the polynomial ring in $N$ variables
representing $\PP^{N-1}$.
In this setting, we show that 
if the free resolution of the ideal $I$ generated by $V$
begins with at 
least $\lceil n/2\rceil$ linear steps, then
the annihilator of $\coker\psi$ is the ideal of 
forms in $\PP(V)$ that vanish on $\alpha_V(\PP^{n-1})$.

If $I$ is an ideal generated in degree $d$, and
$I^k=\mm^{kd}$, then the number of generators $\mu$
of $I$ must satisfy 
$$
{\mu+k-1\choose k} \geq {n+kd-1\choose n-1}.
$$
By \ref{half-way linear2},
this relation is satisfied with $k=2$ if
the resolution of $I$ is linear for $\lceil(n-1)/2\rceil$
steps, and \ref{Us-reg} implies
further
lower bounds. In \ref{almost lin}
we give a stronger lower bound for the number of 
generators of an ideal whose resolution is linear
for $n-2$ steps (the ``almost linear'' case.) It
might be interesting to interpolate to other cases as well.

\medskip
\noindent{\bf The truncation principle}
\smallskip

\noindent 
Since the focus of this paper is on linearly presented ideals,
we have stated many results only for this case. However, 
it is possible to make any ideal $I$ into an ideal with
linear resolution for $s$ steps by truncating, and thus
generalize many of the results. Rather than doing this
throughout the paper, we illustrate it here. The following
result is elementary:
\proposition{truncation} The ideal $J= I\cap \mm^{t_s(I)-s}$
has linear resolution for $s$ steps, while for $p\geq s$ we
have $t_p(J)=t_p(I)$. \Box

For example, \ref{truncation} allows us to deduce \ref{Us-reg} $(a)$
from \ref{Us}.

\medskip
 We are indebted 
to Giulio Caviglia, Aldo Conca, J\"urgen Herzog, 
and Frank Schreyer
for
helpful conversations.

\section{deg syz} Degrees of syzygies

Throughout this paper, $K$ is a field and $S=K[x_1,\dots, x_n]$ is a 
polynomial ring in $n$ variables,
graded with $\deg x_i=1$ (but see remark \ref{nonstandard grading} 
below for the case
of general grading.) We write $\mm=(x_1,\dots,x_n)$ for the
homogeneous maximal ideal of $S$. All tensor products and Tor modules
are taken over the ring $S$.
The Krull dimension of a module $A$ is denoted $\dim A$  (we
use $\dim_K$ for vector space dimension.) 

We write $\reg A$ for the (Castelnuovo-Mumford) regularity of 
a graded $S$-module $A$ (see for example Eisenbud [2004]).
If $A$ is a finitely generated graded vector
space, or more generally an Artinian graded $S$-module,
then $\reg A= \sup \{i \mid A_i\neq 0\}$. If $A$ is a finitely
generated graded $S$-module then
$\reg A$ is defined in terms of local cohomology by the 
formula
$$
\reg A = \max_j \, \{\reg \H^j_\mm (A)+j\}.
$$
For example, if $A=0$ then $\reg A=-\infty$.
We may also compute
$\reg A$ in terms of Tor (or in terms of a minimal 
free resolution) by the formula
$$
\reg A = \max_k\, \{t_k(A) -k\}.
$$
{}From local duality one see that the two ways of expressing the regularity
are also connected ``termwise'' by the inequality
$ t_k(A) -k\geq \reg \H^{n-k}_\mm (A)+n-k$.

The numbers $\reg \H_\mm^j(A)+j$
and  $t_k(A) -k$ will appear often in our
formulas. The next two theorems express the basic technical
result of this paper. 

\theorem{tor-bound1}
Suppose that $A$ and $B$ are finitely generated graded $S$-modules
such that 
$\dim \Tor_1(A,B)\leq 1$, 
and let $j,k$ be integers. If $p\leq \codim A$, $q\leq \codim B$
and $p+q=n-j+k$ then
$$
\reg \H^j_\mm (\Tor_{k}(A,B)) 
\leq 
 t_{p}(A)+  t_q(B)-n.
$$
\smallskip

\theorem{generalization of regularity}
Suppose that $A$ and $B$ are finitely generated graded $S$-modules
such that $\dim \Tor_1(A,B)\leq 1$. 
If $n-j+k\geq \codim A+\codim B$, then
$$
\reg \H^j_\mm (\Tor_{k}(A,B))  
\leq \max_{p+q=n+k-j\atop {p\geq \codim A\atop q\geq \codim B}}
\biggl\{  t_{p}(A) +  t_q(B)\biggr\}-n.
$$
\smallskip

In fact, both these theorems follow from a more general 
statement:

\theorem {tor-bound}
Suppose that $A$ and $B$ are finitely generated graded $S$-modules
such that 
$\dim \Tor_1(A,B)\leq 1$, 
and let $j,k$ be integers. 
For any integers $p,q$ with $p+q=n-j+k$
$$
\reg \H^j_\mm (\Tor_{k}(A,B)) 
\leq \max\{X,Y,Z\}
$$
where
$$\eqalign{
X &= t_{p}(A) +  t_q(B)-n,\cr
Y &=\max_{p'+q'=n-j+k\atop {p'>p}}
\biggl\{ t_{p'}(A) + \reg \H^{n-q'}_\mm(B)\biggr\},\cr
Z &=\max_{p'+q'=n-j+k\atop {p'<p}}
\biggl\{\reg \H^{n-p'}_\mm(A)+ t_{q'}(B)\biggr\}.
}
$$

\noindent{\it Proof of \ref{tor-bound1}.\/}
Since $n-j+k\leq \codim A+\codim B$,  
$q'<\codim B$ in the expression for $Y$
and $p'<\codim A$ in the expression for $Z$,
 so the local cohomology modules in the
expressions for $Y$ and $Z$ in 
\ref{tor-bound} are zero. 
Because
the regularity of the module 0 is
$-\infty$ we have $Y=Z=-\infty$, and
\ref{tor-bound} reduces to \ref{tor-bound1}.\Box\smallskip

\noindent{\it Proof of \ref{generalization of regularity}.\/}
Since $n-j+k\geq \codim A+\codim B$,
we can pick
$p,q$ with  $p\geq \codim A,\ q\geq \codim B$
and $p+q=n-j+k$. 
Replacing
the terms $\reg \H^{n-q'}_\mm(B)+n-q'$ in $Y$ with the possibly larger
terms $ t_{q'}(B)-q'$ (and similarly for $Z$)
in \ref{tor-bound}, we obtain \ref{generalization of regularity}.
\Box

We postpone the proof of \ref{tor-bound} to later in this section.

\remark{nonstandard grading} These formulas adapt easily to the case
where the degrees of the $x_i$ are not assumed to be 1: Setting
$\sigma=\sum\deg x_i$ we must
add $n-\sigma$ to the term $X$ in the main
theorem, and we correspondingly add
$n-\sigma$ to the right hand side of
the formulas in \ref{tor-bound1} and \ref{generalization of
regularity}. The proofs use the comparison
$ t_k(A)-k\geq \reg \H^{n-k}_\mm(A)+\sigma-k$.

Finally, in case the module $B$ is Cohen-Macaulay, a special case
of the inequality takes on a simple form no matter what the 
relation of $n+j-k$ and $\codim A+\codim B$:

\corollary{CM case} Suppose that $A$ and 
$B$ are finitely generated graded $S$-modules
such that $\dim \Tor_1(A,B)\leq 1$. If $B$ is
Cohen-Macaulay of dimension $b$ then
$$
\reg \H^j_\mm (\Tor_{k}(A,B)) 
\leq 
 t_{b-j+k}(A)-b + \reg B.
$$

\noindent For example, when $B$ has finite length,
this statement reduces to the easy formula
$\reg(\Tor_k(A, B))\leq t_k(A)+\reg B$.
\smallskip

\noindent{\sl Proof of \ref{CM case}} 
Take $q=n-b=\codim B$ in \ref{tor-bound}.
Because $B$ is Cohen-Macaulay, the only
nonvanishing local cohomology of $B$ is 
$\H^b_\mm(B)$. 
The terms $\reg \H^{n-q'}$ 
that appear in the expression for $Y$ in \ref{tor-bound}
are all $-\infty$ because when $p'>p$ we have $n-q'>b$,
so $\H^{n-q'}_\mm(B)=0.$ The terms
$ t_{q'}(B)-q'$ that appear in the
expression for $Z$ are all $-\infty$
because when $p'<p$ the number $q'$ is bigger than
$n-b$, the projective dimension of $B$.
\Box
\smallskip

\ref{tor-bound} can fail without the assumption that
$\dim\Tor_1(A,B)\leq 1$,
even in the case where $A=B=R/I$ is 2-dimensional and $n=4$:
for instance it is easy to check that
the example of Conca given after \ref{reg of powers}
(with $r=2$, say) does not satisfy \ref{reg of Tor}.
The assumption 
$\dim \Tor_1(A,B)\leq 1$ is used in the proof to ensure the
degeneration of a certain spectral sequence. In fact, to
achieve the vanishing needed, it is enough to make the
weaker assumption that 
$$\eqalign{
\H^{j+k+1}_\mm(\Tor_{i+k}(A,B))&=0\cr
\H^{j-k-1}_\mm(\Tor_{i-k}(A,B))&=0
}$$
for all $k\geq 1$. This ensures that all the differentials
in the $E_k$ pages of the spectral sequence in the proof
that come from or go to this term are zero, so the term, and
not
a proper subquotient of it, is a subquotient of the 
corresponding term in the limit of the spectral sequence.

We note that the hypothesis $\dim \Tor_1(A,B)$ is always
satisfied if $A,B$ are dimensionally transverse in the
sense that $\codim A\otimes B \geq \codim A+\codim B$
(in which case equality holds) and $A,B$ are both
locally Cohen-Macaulay off a set of codimension $\geq 2$.
\smallskip\goodbreak

For any graded $S$-module we write 
$\mindeg T= \inf \{i \mid T_i\neq 0\}$. If 
$T=0$ we set $\mindeg T=\infty.$
\smallskip

\noindent{\sl Proof of \ref{tor-bound}}
Let $\FF:\cdots \to F_1\to F_0$ be a minimal free resolution of $A$ and let
$\GG:\cdots \to G_1\to G_0$ be a minimal free resolution of $B$. The proof
consists of an analysis of the double complex
$\FF^*\otimes \GG^*=(\FF\otimes \GG)^*$
where $*$ denotes $\Hom(-,S)$.

For any finite complex $\KK: \cdots \to K_n\to K_{n-1}\to\cdots$ 
of free $S$-modules there is a spectral sequence
with $E_2$ term
$
\Ext_S^s(\H_t(\KK),S)
$
converging to $\H^{s+t}(\KK^*)$, obtained
from the double complex $\Hom(\KK, \II)$, where
$\II$ is an injective resolution of $S$. 
We apply this to $\KK={\rm Tot}(\FF\otimes \GG)$.
Since $\Tor_1(A,B)$ has Krull dimension at most 1,
Auslander's Theorem [1961] on the rigidity of
Tor
shows that 
$\H_t(\FF\otimes\GG)=\Tor_t(A,B)$ has dimension $\leq 1$
for every $t\geq 1$. It follows that
$
\Ext_S^s(\H_t(\KK),S)
$
is nonzero only when $t=0$ or when $s=n-1$ or $s=n$. The $E_2$
differential $\Ext_S^s(\H_t(\KK),S)\to \Ext_S^{s+2}(\H_{t-1}(\KK),S)$
thus vanishes and the spectral sequence degenerates
at $E_2$. The degeneracy in turn shows that
$\Ext_S^s(\H_t(\KK),S)$
is a subquotient of  $\H^{s+t}(\KK^*)$.

By local duality
$$\eqalign{
\H^j_\mm(\Tor_k(A,B)) &= \H^j_\mm(\H_k(\KK))\cr
&=\Hom_K(\Ext^{n-j}(\H_k(\KK),S),K)(n)
}
$$
where $\Hom_K$ denotes the set of graded homomorphisms.
Since $\Ext^{n-j}(\H_k(\KK),S)$ is a subquotient of 
$\H^{n-j+k}(\KK^*)$, it follows that
$$
\reg \H^j_\mm \Tor_{k}(A,B)
\leq
-\mindeg \H^{n-j+k}(\KK^*)-n.
$$

To prove  \ref{tor-bound} we need to show that 
any homogeneous element
$\zeta\in \H^{n-j+k}(\KK^*)$ of degree
$$
\deg \zeta< -\max\{X,Y,Z\}-p-q
=-\min\{p+q-X,p+q-Y,p+q-Z\}
$$ 
is zero.
We have
$$
p+q-X=- t_p(A)- t_q(B)
$$
and by local duality
$$
\eqalign{
{(*)}\qquad p+q-Y &= \min_{p'+q'=n-j+k\atop {p'>p}}
\biggl\{- t_{p'}(A) + \mindeg \Ext^{q'}(B,S)
\biggr\},\cr
{(**)}\qquad
p+q-Z &=
\min_{p'+q'=n-j+k\atop {p'<p}}
\biggl\{\mindeg \Ext^{p'}(A,S)- t_{q'}(B)\biggr\}.
}$$
Let $z=\{z^{p',q'}\mid p'+q'=p+q\}$ be a homogeneous cycle 
of $\KK^*$ representing
$\zeta$.
Since 
$$\eqalign{
\mindeg (F_p^*\otimes G_q^*) &=
\mindeg (F_p^*\otimes K)+\mindeg(G_q^*\otimes K)\cr
&=
- t_p(A) - t_q(B)\cr
&>\deg \zeta,
}$$
it follows that $z^{p,q}=0$.
To finish the proof we will show that the
other components $z^{p',q'}$ are also zero.

By equation (**) the vertical homology of $\KK^*$ at
$(\KK^*)^{p',q'}$ is zero in degree $\deg \zeta$
when 
$p'+q'=p+q$ and $p'<p$, while by equation (*) the horizontal
homology of $\KK^*$ is zero at $(\KK^*)^{p',q'}$
in degree $\deg \zeta$ when
$p'+q'=p+q$ and $p'>p$.

We may thus complete the proof by applying the
following 
more general Lemma to the complex $\LL$ formed
by taking the degree $\deg \zeta$ part of $\KK^*$.
The result gives information about the total
cycles in the double complex
$$\LL:\qquad\qquad
\diagram 
&\ddots&\uTo_{d_{\rm vert}}&&\uTo_{d_{\rm vert}}
\cr
&\rTo^{d_{\rm hor}}&L^{p',q'}&\rTo^{d_{\rm hor}}&
L^{p',q'+1}&\rTo^{d_{\rm hor}}&
\cr
&&\uTo_{d_{\rm vert}}&\ddots&\uTo_{d_{\rm vert}}
\cr
&\rTo^{d_{\rm hor}} & L^{p'-1,q'}&\rTo^{d_{\rm hor}}&
L^{p'-1,q'+1}&\rTo^{d_{\rm hor}}&
\cr
&&\uTo_{d_{\rm vert}}&&\uTo_{d_{\rm vert}}&\ddots
\enddiagram
$$
\smallskip

\lemma{sliding} Let $\LL$ be any bounded below double complex,
with noation as above,
suppose that $p,q$ are chosen so that the vertical
homology of $\LL$ is zero at $L^{p',q'}$ when
$p'+q'=p+q$ and $p'<p$, and the horizontal
homology of $\LL$ is zero at $L^{p',q'}$ when
$p'+q'=p+q$ and $p'>p$. If 
$\zeta \in H^{p+q}Tot(\LL)$ represented by a cycle
$$
z=(z^{p',q'})\in \oplus_{p'+q'=p+q}L^{p',q'}
$$
satisfies $z^{p,q}=0$, then $\zeta=0$.

\proof 
We have $d_{\rm vert}(z^{p-1,q+1})= -d_{\rm hor}z^{p,q} = 0$.
By our assumption the vertical homology vanishes at $L^{p-1,q+1}$
so $z^{p-1,q+1}=  d_{\rm vert}(w)$ for some $w\in L^{p-2,q+1}$.
Subtracting $d_{\rm Tot} w $  from $z$ we get a homologous cycle
$y$ whose components $y^{p',q'}$ agree with $z^{p',q'}$ for
$p'\geq p$, but $y^{p-1,q+1}=0$. Repeating this process we see
that $z$ is homologous to a cycle $x$ with $x^{p',q'}=z^{p',q'}$
for $p'\geq p$ while $x^{p',q'}=0$ for $p'<p$.

Similarly, using the fact that the horizontal homology
is zero at $L^{p',q'}$ for $p'>p$ and $p'+q'=p+q$,
we can change $x$ by a boundary to arrive at a cycle that is 0
in every component, so $\zeta=0$.
\Box
\smallskip

In the special case where $B$ is a Gorenstein
factor ring of $S$
we can describe when 
\ref{tor-bound} (in the form of \ref{CM case}) is sharp.
Suppose $\phi:F'\to F$ is a map of 
graded free modules such that $\reg F = d$. By a 
\emph{generalized row of $\phi$ of maximal degree} we mean
the composition of $\phi$ with a projection $F\to S(-d)$.
By ``the entries'' of this row we mean the ideal that is
the image of the corresponding map $F(d)\to S$.

\proposition{sharpness}
Suppose that $A$ is a finitely generated graded $S$-module
with free resolution 
$$
\cdots \rTo F_t\rTo^{\phi_t} F_{t-1}\rTo \cdots \rTo^{\phi_1} F_0
$$
and  $J$ is an ideal such that $S/J$ is Gorenstein of dimension $b$
and 
$A/JA$ has finite length. If $k\leq \codim A -b$ then
$$
\reg \Tor_k(A, S/J) \leq  t_{b+k}(A)-b+ \reg S/J
$$
with equality  if and only if $J$ contains the ideal
generated by the entries in 
some generalized row of maximal degree of
$\phi_{b+k+1}$.

\proof
The inequality is \ref{generalization of regularity}.
Since $B=S/J$ is Cohen-Macaulay we have
$\reg B= t_{n-b}(B)-n+b$. Since $A\otimes B=A/JA$
has finite length,
$$
\reg \Tor_k(A,B) = -\mindeg\Hom_K(\Tor_k(A,B), K).
$$
By local duality, we can rewrite this as 
$-(\mindeg \Ext^n(\Tor_k(A,B),S))-n.$

We now use the 
notation and spectral sequence from the proof of 
\ref{tor-bound}. 
Because $A\otimes B$ has finite length,
the $E_2$ page of the spectral sequence for the
homology of $\KK^*$ has
nonzero terms only in one row and one column, and
if follows that 
$\Ext^n(\Tor_k(A,B),S)=\H^{n+k}(\KK^*).$

{}From this we see that equality holds in
\ref{sharpness} if and only if 
$\mindeg \H^{n+k}\Tot (\FF^*\otimes \GG^*) = 
\mindeg(F_{b+k}^*\otimes G_{n-b}^*)$.
Because $B$ is Gorenstein we may write $G_{n-b}^*=S(e)$ 
for some $e$.
Moreover $\GG^*$ is a resolution of 
$\Ext^{n-b}(B,S) = B(e)$.
It follows that
$\H^{n+k}\Tot (\FF^*\otimes \GG^*) 
\iso \H^{b+k}(\FF^*\otimes B)(e)$.
Hence equality holds if and only if 
$\mindeg(F_{b+k}^*\otimes S(e)) 
= \mindeg \H^{b+k}(\FF^*\otimes B)(e)$.
Since $\FF^*$ is a minimal complex, this is equivalent to saying that
a generator of minimal degree of $F_{b+k}^*$ is a cycle mod $J$; 
that is, 
$J$ contains the ideal
generated by the entries in some generalized row of maximal degree of
$\phi_{b+k+1}$.
\Box

\section{C-M reg section} Castelnuovo-Mumford Regularity

The following is an extension of results of Sidman [2002]
and Caviglia [2004], who treat the case $k=0$ by different
methods. 

\corollary{reg of Tor} If $A$ and $B$ are finitely generated
 graded $S$-modules
such that $\dim \Tor_1(A,B) \leq 1$, then
$$
\reg \Tor_k(A,B) \leq \reg A +\reg B +k.
$$

\proof We use the formula
$$
\reg M = \max_j\{\reg \H^j_\mm(M) +j\mid j\geq 0\}
$$
to compute $\reg \Tor_k(A,B)$, and 
$$
\reg A+\reg B=
\max\{ t_p(A) -p+  t_q(A) -q\mid p,q\geq 0\}.
$$ 
The proof
is then  a straightforward application of the
inequalities in Theorems \refn{tor-bound1} and
\refn{generalization of regularity}.
\Box

\corollary{newregtor} Suppose that $A$ and $B$ are finitely generated
graded $S$-modules
such that $\dim \Tor_1(A,B) \leq 1$, and let $k$ be an integer.
If $k+ \dim B \leq p \leq \codim A$ then
$$
\reg \Tor_k(A,B) \leq  t_p(A) +  t_{n+k-p}(B) -n.
$$

\proof
Since $p\leq \codim A$ and $n+k-p\leq \codim B$, 
\ref{tor-bound1} gives
$$
\reg\Tor_k(A,B)\leq \max_{j=0,1}\{t_p(A)+t_{n-j+k-p}(B)+j-n\}.
$$
But $t_{n-j+k-p}(B)+j \leq t_{n+k-p}(B)$, again because
$n+k-p\leq \codim B$. 
\Box

\corollary{reg of Tor Cohen-Macaulay} Suppose that
 $A$ and $B$ are graded $S$-modules
such that $\delta:=\dim \Tor_1(A,B) \leq 1$. If  
$B$ is a Cohen-Macaulay module of dimension $b$, then
for $k>0$
$$\eqalign{
\reg \Tor_k&(A, B)\cr
\leq &
\max\{ t_p(A) - p\mid b+k-\delta\leq p\leq b+k\}\cr
&+\reg B +k.
}$$

\proof 
Notice that $\dim \Tor_k(A,B)\leq \delta$ by the rigidity
of Tor (see Auslander [1961]). Thus the assertion follows
from \ref{CM case}.
\Box

As an application
of Corollaries \refn{reg of Tor} and \refn{reg of Tor Cohen-Macaulay} 
with $k=1$, we have

\corollary{intersection and product}
If $I$ and $J$ are homogeneous ideals of $S$ such that
$(IJ)_d = (I\cap J)_d$ for $d>>0$, then the equality holds
for all
$d\geq \reg I +\reg J$. If in addition $S/J$ is Cohen-Macaulay
of dimension $b$, then it suffices that
$$
d \geq  t_b(I)-b + \reg J.
$$

\proof We use the formula
$\Tor_1(S/I, S/J)=(I\cap J)/IJ$, and
apply Corollaries \refn{reg of Tor} and \refn{reg of Tor Cohen-Macaulay}.
\Box

Suppose that $X, Y\subset \PP^{n-1}$ are schemes. The ideal
$I_{X\cap Y}$
of $X\cap Y$ is the saturation of the sum of the ideals of $X$ and
$Y$; that is, they agree in high degrees.
 Using Theorems \refn{tor-bound1} and \refn{generalization of regularity}
we can make this quantitative in the case where $X$ and $Y$ meet
at most in dimension 0. Note that in this
case $\codim X+\codim Y\geq n-1$.

\corollary{sum of ideals} 
Let $X,Y\subset \PP^{n-1}$
be schemes with ideals $I,J\subset S$. Suppose that $\dim X\cap Y = 0$.
 \item{$(a)$} If $\codim X +\codim Y\geq n$, then any form of
degree $d$ vanishing on $X\cap Y$ is a sum of a form
vanishing on $X$ and a form vanishing on $Y$ as long as
$$
d>  t_p(S/I) +  t_q(S/J) -n 
$$
for some integers $p,q$ satisfying
$p\leq \codim X,\ q\leq \codim Y,$ and $p+q=n.$
\item{$(b)$} If $\codim X +\codim Y= n-1$, then any form of
degree $d$ vanishing on $X\cap Y$ is a sum of a form
vanishing on $X$ and a form vanishing on $Y$ as long as
$$\eqalign{
d> \max \{ 
& t_{1+\codim X}(S/I) +  t_{\codim Y}(S/J),
\cr 
& t_{\codim X}(S/I) +  t_{1+\codim Y}(S/J)
\}-n.
}
$$

\proof Notice that
$S/(I+J)=(S/I)\otimes (S/J)$. It follows that $S/(I+J)$
is saturated in degree $d$ if 
$\H^0_\mm(\Tor_0(S/I, S/J))_d=0$. Cases (a) and (b)
follow from Theorems \refn{tor-bound1} 
and \refn{generalization of regularity},
with $j=k=0$.
\Box

A similar result follows for any schemes $X$ and $Y$ 
whose intersection is ``homologically transverse''
except along a zero-dimensional set in $\PP^{n-1}$
(but the sum of the
codimensions of $X$ and $Y$ may then be $<n-1$, in which
case more terms appear in case (b)).

\section{subadd section} Convexity of degrees of syzygies

\ref{tor-bound1}
yields a kind of ``triangle inequality"
or convexity for degrees of syzygies
that seems to be new even
in the case where $A=B$ is a module of finite length.

\corollary{subadd}
Suppose that $A$ and $B$ are finitely generated graded $S$-modules
such that $\dim \Tor_1(A,B)\leq 1$, then 
$$
t_n(A\otimes B) \leq  t_{p}(A) +  t_{n-p}(B)
$$
whenever $\dim B\leq p \leq \codim A$. In particular, if $A=B=S/I$
is a cyclic module of dimension $\leq 1$,
then the function $p\mapsto t_p(S/I)$ satisfies the 
weak convexity condition
$$
 t_n(S/I) \leq  t_{p}(S/I) +  t_{n-p}(S/I).
$$
for $0\leq p\leq n$.

When $\dim B>\codim A$ a similar result follows from 
\ref{generalization of regularity}.

\proof For any finitely generated graded module $M$, 
$$
\Tor_n(M,K) = \ker\biggl( M(-n)\rTo^{\pmatrix {x_1\cr\vdots\cr x_n}} 
M^n(-n+1)\biggr)
= {\rm socle }\, M(-n).
$$
as can be calculated from the Koszul resolution of $K$.
Thus $\reg \Tor_n(A\otimes B)= \reg \H^0_\mm(A\otimes B)+n$,
and the assertion follows from \ref{tor-bound1}.
\Box

If a module $A$ is annihilated by 
an $\mm$-primary ideal $J$, then it is immediate that
the degree of the socle of $A$ is bounded above by the
highest degree of a generator of $A$ plus the highest degree of the 
socle of $S/J$. This relation can be written as
$t_n(A)\leq t_0(A) + t_n(S/J)$.
The following result gives such a bound without the assumption
that $J$ is $\mm$-primary.

\corollary{socle estimation} Suppose that $A$ is a
finitely generated graded
$S$-module of codimension $c$ and that $\delta:=\dim A-\depth A\leq 1$.
Let $J$ be a homogeneous ideal contained
in the annihilator of $A$. If 
$\depth S/J\geq \depth A$ 
then for $0\leq q\leq \codim J$
$$
 t_{c+\delta}(A)\leq  t_{c+\delta-q}(A) + t_q(S/J).
$$
In particular:
\item{$($a$)$} If the annihilator of $A$ contains a 
regular sequence of forms of degree
$d_1\leq\cdots\leq d_q$  then
$$
 t_{c+\delta}(A) \leq  t_{c+\delta-q}(A)+d_1+\cdots +d_q.
$$
\item{$($b$)$} If $J$ is perfect and is generated in degree $d$
with linear resolution, then 
$$
 t_{c+\delta}(A) \leq  t_{c+\delta-q}(A)+d+q-1.
$$

\proof We may harmlessly assume that $K$ is
infinite. If $\dim A>1$ a general sequence of
$\depth A$ linear forms is a regular sequence
on both $A$ and $S/J$, so we factor out these
linear forms (and work over the corresponding 
factor ring of $S$) without changing the statement.
Thus we may suppose $\dim A\leq 1$ and $\depth A=0$,
so $n=c+\delta$. 
Since the case $q=0$
is trivial, we may suppose that $q\geq 1$.

We now apply 
\ref{tor-bound1} 
with $k=j=0$, $B=S/J$ and $p=n-q$. As $p\leq \codim A$
 we obtain
$\reg \H^0_\mm(A)\leq  t_p(A)+ t_q(B)-n$. 
Since
$\reg \H^0_\mm(\Tor_0(A,S/J))=\reg \H^0_\mm (A) =
 t_n(A)-n$, this gives
the first statement. Parts (a) and (b) follow immediately
by computing $\Tor(S/J,K)$ in the given cases.
\Box
 
\example{} If $X$ is an arithmetically
Cohen-Macaulay scheme of codimension $c$ in $\PP^m$, 
with ideal $I$ and $X$ is contained in a nondegenerate variety
of codimension $q$ and (minimal) degree $q+1$,
then by part (b) of \ref{socle estimation},
$$
 t_{c}(S/I)\leq  t_{c-q}(S/I)+q+1.
$$
\smallskip
\example{caviglia} (G.~Caviglia, Thesis) The principle of part (a) 
of \ref{socle estimation} does not
hold for individual steps in the resolution. For example, if
$$
I=(x_1^3,\dots, x_4^3, (x_1+\cdots +x_4)^3)\subset S=K[x_1,\dots, x_4]
$$
then $ t_1(S/I)=3$ while $ t_2(S/I)=7 > 3+3$.

Similarly, if 
$$
I = (x_1^n, x_2^n, x_1x_3^{n-1}-x_2x_4^{n-1})\subset S=K[x_1,\dots, x_4],
$$
then $ t_1(S/I)=n$ 
while $ t_2(S/I) = n^2 > 2n$ for $n\geq 3$, and in fact
$\reg(S/I)=n^2-2$.

The same idea shows that the dimension bound on $\Tor_1(A,B)$ 
is necessary in
\ref{CM case} and \ref{reg of Tor}. 
In the ring $T = S[t]$, we can write $(I,t) = J + L$ where
$J = (x_1^n, x_2^n, x_1x_3^{n-1}-x_2x_4^{n-1}+t^n)$ and $L = (t)$. Note
that both $J$ and $L$ are complete intersections. For $n\geq 3$, 
$$
\reg \Tor_0(T/J, T/L) = \reg (S/I) = n^2-2 > \reg (T/J) + \reg (T/L) = 3n-3.
$$
In this case $\dim \Tor_1(T/J, T/L) = 2$.

\section{specializing section} Specialization and degrees of syzygies

As an application of \ref{reg of Tor Cohen-Macaulay}
we give a bound for the saturation and regularity of a plane section:
If we take the
case where $I$ is the saturated ideal of $X$, and $Y$ is a linear
space, we obtain a result that generalizes Theorem 1.2
of Eisenbud-Green-Hulek-Popescu [2004a].

\corollary{regularity of plane sections}
Let $X\subset \PP^{n-1}$ be a scheme, and let $\Lambda\subset \PP^{n-1}$
be a linear subspace such that the sheaf $\Tor_1(\O_X, \O_\Lambda)$ is 
supported on a finite set. 
Let $I\subset S$ be any homogeneous ideal defining $X$,
and let $L\subset S$ be the ideal of $\Lambda.$
\item{$($a$)$} 
The restriction map
$$
I_d= \H^0(\I_X(d))\to \H^0(\I_{X\cap\Lambda, \Lambda}(d))
$$
is surjective for all 
$d\geq  t_{\dim \Lambda}(I)-\dim\Lambda$.
\item{$($b$)$}
Let $c$ be the codimension of $X\cap\Lambda$
in $\Lambda$. We have
$$\eqalign{
\reg (\I_{X\cap \Lambda})&= \reg ({\I_X+I_\Lambda\over \I_\Lambda})\cr
&\leq 
\max
\{ t_{p}(I) -p\mid c-1\leq p\leq \dim \Lambda -1\}.
}
$$

The hypothesis that the sheaf $\Tor_1(\O_X, \O_\Lambda)$ is supported
on a finite set is satisfied for general $\Lambda$ of
any dimension, or for any $\Lambda$ such that $X\cap\Lambda$ is finite.

\proof 
By \ref{CM case} we have 
$$\eqalign{
\reg\H^j_\mm &(S/(I+L))\cr
&=\reg \H^j_\mm (\Tor_0(S/I, S/L))\cr
&\leq t_{\dim(S/L) -j}(S/I)-\dim S/L\cr 
&= t_{\dim \Lambda-j}(I)-\dim\Lambda-1\cr
&< t_{\dim \Lambda-j}(I)-\dim\Lambda.}
$$
Taking $j=0$ in the inequalities,
we see that  $I+L$ is saturated in 
degree $d$ when $d\geq  t_{\dim \Lambda}(I)-\dim\Lambda$,
proving part (a).
Adding $j$ to both sides and taking the maximum over $j$ 
for $1\leq j\leq \dim S/(I+L) = \dim X\cap\Lambda +1$
we see that 
$$\eqalign{
\reg\ & \I_{ X\cap \Lambda} \cr
&= \max_{1\leq j}\{ \reg \oplus_m \H^j (\I_{X\cap\Lambda}(m))+j+1\} \cr
&=\max_{1\leq j\leq \dim S/(I+L) } \{\reg \H^j_\mm (S/(I+L))+j+1\}\cr
&\leq \max_{1\leq j\leq \dim S/(I+L) } 
\{ t_{\dim \Lambda-j}(I)-\dim\Lambda+j\},
}
$$
which is the desired inequality.
\Box

We say that the resolution of a finitely generated
graded $S$-module $A$ generated in
a single degree $d$ is \emph{linear for $q$ steps} if it has the
form
$$
\cdots \rTo S^{n_q}(-d-q)\rTo\cdots\rTo S^{n_0}(-d)\rTo A\rTo 0.
$$

\corollary{linear complements} 
Suppose that $I\subset S$ is a homogeneous ideal, let $p$ be an 
integer and set
$m= t_p(S/I)$.
Let $L \subset S$ be any ideal generated by $n-p$ independent linear forms.
If $I+L$ contains a power of $\mm$ (which will
always be true if $K$ is infinite, $L$ is
general and $p \leq \codim I$) then $I+L$
contains $\mm^{m-p+1}$, and more generally 
$$
\mm^{m-p+s}\subset I+L^s.
$$
For example, if $I$ is generated in degree $d$ and 
the minimal free resolution of $I$ is linear for
$p-1$ steps, then 
$$
\mm^d\subset I+L.
$$

\proof The resolution of $L^s$ is linear, 
as one can see by computing the degree of the socle
of $S/L^s$ (in fact, the resolution can be computed as
an Eagon-Northcott complex, see Eisenbud [1995], pg. 600). Hence
$ t_{n-p}(S/L^s)= n-p+s-1$. As $p \leq \codim I$, 
\ref{tor-bound1} gives $\reg \H^0_\mm(S/I \otimes S/L^s) 
\leq m-p+s-1$, which is the asserted result.
\Box

Notice that the containment $\mm^d\subset I+L$ in \ref{linear complements}
actually gives that $I$ and $\mm^d$ coincide modulo $L$.

\corollary{initials} 
Suppose that $I\subset S$ is a homogeneous $\mm$-primary ideal,
and let $\In I$ denote the initial ideal of $I$ with
respect to the reverse lexicographic order on the monomials of $S$.
If $m= t_p(S/I)$ then
$$
(x_1,\dots,x_{p})^{m-p+1} \subset \In I
$$
In particular,
if $I$ is generated in degree $d$
and the resolution of $I$ is linear for
$p-1$ steps, then the initial ideal of $I$ in reverse lexicographic
order contains $(x_1,\dots,x_{p})^d$.

\proof \ref{linear complements} shows that 
$
\mm^{m-p+1}\subset I+L,
$
where $L = (x_{p+1},\dots,x_n)$. Because the monomial order
is reverse lexicographic, $\In (I+L) = (\In I)+L$ (see Eisenbud [1995],
Proposition 15.12). Thus $\mm^{m-p+1}\subset(\In I)+L$, whence $(x_1,\dots,x_{p})^{m-p+1}
\subset \In I$.
\Box

In the case where $I$ is
$\mm$-primary and linearly
presented, \ref{initials} says that $(x_1,x_2)^d\subset \In I$.
In generic coordinates we hope for a stronger
inclusion:

\conjecture{gin conjecture 1}
Suppose that the ideal $I\subset S$ is
$\mm$-primary, linearly
presented, and generated in degree $d$.
If $K$ is infinite, then
$$
\mm^d\subset I+(z_3,\dots,z_n)^2
$$ 
for sufficiently
general linear forms $z_3,\dots,z_n$,
or equivalently 
$$
(z_1,z_2)^{d-1}\mm\subset \gin I,
$$
where $\gin I$ denotes the reverse lexicographic initial ideal with
respect to generic coordinates $z_1,\dots,z_n$. If the resolution of
$I$ is linear for $p$ steps, then we similarly conjecture that
$$
\mm^d\subset I+(z_{p+2},\dots,z_n)^2
$$ 
for sufficiently
general linear forms $z_i$.

We were led to this conjecture studying \ref{Eisenbud-Ulrich}. 
In case $n=3$ and $S/I$ is
Gorenstein, \ref{gin conjecture 1} follows from the Hard Lefschetz
property proved by
Harima, Migliore, Nagel and Watanabe [2003].
We have observed it experimentally in a large number of
other cases with $n=3$ and $n=4$. 

\corollary{HeHi}
Suppose that $K$ has characteristic zero and $I \subset S$ is a 
homogeneous $\mm$-primary ideal. If $I$ is generated in degree $d$
and the resolution of $I$ is linear for $n-2$ steps, then
$\mu(\gin I)=\mu(\mm^d)$. 

\proof  \ref{linear complements} shows that $I+(z)=\mm^d+(z)$ 
for every linear form $z$ in $S$. But then $\mu(\gin I)=\mu(\mm^d)$
by Herzog and Hibi [2003]. 
\Box

\section{quadrics section} Ideals generated by quadrics

If an ideal $I$ generated in degree
$d$ has a resolution that is linear for
 $q$ steps, then
by \ref{linear complements} we have
$\mm^d\subset I+(\ell_{q+2},\dots,\ell_n)$ for every set of independent
linear forms $\ell_{q+2},\dots,\ell_n$. For ideals generated
by quadrics, this latter condition is easy to interpret.
For simplicity we assume throughout this section that the base
field $K$ is algebraically closed of characteristic not 2. We will
identify a quadric and its associated symmetric bilinear form.

Recall that a $m$-dimensional vector space of quadrics in $n$ variables (with a basis)
can be described by a symmetric $n\times n$ matrix of linear forms
in $m$ variables; to get the symmetric matrix 
corresponding to the $i$-th quadric, just set all but the 
$i$-th variable equal to 0, and set the $i$-th variable equal to 1.
We call a symmetric matrix of linear forms in $m$ variables
{\it symmetrically $q$-generic\/} if every generalized
principal $(q+1)\times (q+1)$ 
submatrix has independent entries on and above the diagonal
(here a principal submatrix is
one involving same rows as columns, and a generalized submatrix
of $A$ is a submatrix of $PAP^*$ for some invertible matrix $P$.)
These definitions are adapted from the notion of $k$-generic matrices
in Eisenbud [1988]. In particular, symmetrically 1-generic matrices
are the same as 1-generic matrices that happen to be symmetric.

It is convenient for our purpose to specify a space of quadrics
via its orthogonal complement. A symmetric matrix $A$ representing a 
quadric may be thought of as a linear transformation $A:W\to W^*$.
The dual of the $\Hom(W,W^*)$ is $\Hom(W^*, W)$ by the pairing
$(A,B)={\rm Trace AB}$. What this means in practice for symmetric
matrices $A=(a_{i,j}),\; B=(b_{i,j})$
is that $(A,B)=\sum_{i,j}a_{i,j}b_{i,j}$. Thus from a space of (quadratic or)
bilinear forms $U$ we can construct a space $U^\perp$ of
(quadratic or) bilinear forms. This is the degree 2 part of the
the ``annihilator ideal'' that appears for example in Eisenbud [1995],
Section 21.2.

The orthogonal complement construction allows
us to give examples of symmetrically
$q$-generic families of quadrics for all $q$:

\proposition{rank and rank} A quadratic form $Q$ 
has rank $\geq q+2$ if and only if the family 
$(Q)^\perp$ of quadratic forms orthogonal
to $Q$ is symmetrically $q$-generic. 

\proof If $Q$ has rank $\leq q+1$ then, after a change of
variables, $Q$ will be represented by a diagonal matrix with at most $q+1$
nonzero entries. 
It follows that
the matrices in $(Q)^\perp$ satisfy a linear equation
among the entries of some $(q+1)\times (q+1)$ principal
submatrix, so the family is not $q+1$-generic.

Conversely, If the family $V=Q^{\perp}$ is not $q$-generic, then there is
a relation on the entries of a $(q+1)\times (q+1)$
generalized principal submatrix. The coefficients of this
relation define a  quadratic form $Q'$
of rank at most $q+1$ so that
$V\subset (Q')^\perp$. Since both sides are codimension 1 in
$S_2$, they are equal, and it follows that $Q'$ and $Q$ generate
the same 1-dimensional subspace. In particular they have the same rank.
\Box

\proposition{genericity of spaces of quadrics} Let $V\subset S_2$
be an $m$-dimensional vector space of quadrics in 
$n$ variables, represented by
a symmetric matrix $A$ of linear forms in $m$ variables. The ideal
$I$ generated by $V$ has the property that 
$\mm^2\subset I+(\ell_{q+2},\dots,\ell_n)$ for every set of independent
linear forms $\ell_{q+2},\dots,\ell_n$
if and only if $A$ is symmetrically $q$-generic.

\proof The space of quadratic forms 
$\overline V\subset (S/(\ell_{q+2},\dots,\ell_n))_2$
corresponds to the $(q+1)\times (q+1)$
 generalized submatrix of $A$ obtained by leaving
out rows and columns corresponding to the linear forms $\ell_i$.
Its ${q+2\choose 2}$ entries on and above the diagonal are linearly independent if
and only if it corresponds to a space of quadrics of dimension ${q+2\choose 2}$,
which is the dimension of $(S/(\ell_{q+2},\dots,\ell_n))_2$.
\Box

\corollary{bound on number of quadrics} If the ideal $I$ 
generated by $m$ quadratic forms in $n$ variables
is $\mm$-primary and satisfies 
$\mm^2\subset I+(\ell_{3},\dots,\ell_n)$ for every set of independent
linear forms $\ell_{3},\dots,\ell_n$, then $m\geq 2n-1$.

\proof The entries of a 1-generic $n\times n$ matrix must span
a space of at least dimension $2n-1$; see Eisenbud [1988].\Box

\noindent{\bf Example.} The ``catalecticant" (or Hankel) matrix
$$
\pmatrix{
x_1 & x_2 & x_3&\dots\cr
x_2& x_3& \dots\cr
x_3&\dots\cr
\vdots}
$$
is a symmetrically 1-generic matrix representing a
$2n-1$ dimensional space of quadrics.

\corollary{bound on linear length for quadrics} Let $V\subset S_2$
be an $m$-dimensional vector space of quadrics in 
$n$ variables, represented by
a symmetric matrix $A$ of linear forms in $m$ variables. 
If $A$ is not symmetrically $q$-generic, then the ideal
$I$ generated by $V$ has a free resolution with at most $q-1$ linear
steps.\Box

In case $V$ has codimension 1 in the space of all quadrics, 
\ref{bound on linear length for quadrics} is sharp:

\proposition{resolution of quadrics} Let $V\subset S_2$ be 
a codimension 1 subspace of the quadratic forms of $S$.
The ideal generated by $V$
has $q\geq 0$ linear steps in its resolution if and only if
$V$ is $q$-generic. 

\proof Let $Q$ be a quadratic form generating the
 orthogonal complement of $V$. Suppose that the rank
of $Q$ is $q+2$. By \ref{rank and rank} and
\ref{bound on linear length for quadrics},
it suffices to show that the resolution
of $I=(V)$ has  $q$ linear steps.

Let $J$ be the annihilator of $Q$ in the sense of 
Eisenbud [1995], Section 21.2. Thus $S/J$ is Gorenstein
with ``dual socle generator $Q$'', and $J$ contains 
$n-q-2$ independent linear forms $\ell_{q+3},\dots,\ell_{n}$.

If $q+2=n$, the resolution
of $S/J$ has the form
$$
0\to S(-n-2) \rTo \oplus S(-n)\rTo\cdots\rTo \oplus S(-2)\rTo S,
$$
showing that $J=I$ and proving the Proposition in this case.

For arbitrary $q$ we see that the resolution of $S/J$
is the tensor product of a Koszul complex on $n-q-2$ linear
forms with a resolution of $S/J'$ where $S/J'$
is Gorenstein of codimension $q+2$ and has resolution 
similar to the one above. Thus the regularity of $S/J$ is 2, and
$$
\dim \Tor_{n-t}(S/J,K)_{n-t+2} = {n-q-2\choose n-t},
$$
which vanishes for $t> q+2$.

In particular, $J$ is generated in
degrees 1 and 2, so $I$ may be written as $I=J\cap \mm^2$.
We thus get an exact sequence $0\to I\to J\to K(-1)^{n-q-2}$.
Computing $\Tor_{n-t}(S/I,K)$ from this exact sequence,
we see that $S/I$ has  $q$ linear steps in its resolution
as required.
\Box

Using the theory of matrix pencils, it should be possible to analyze all
the complements of codimension two sets of 
quadrics.

\section{reg prods section} Regularity of products and powers

In this section we give our results on \ref{Us}.
At present
we cannot even prove that  $I^2$ has linear
presentation! But we can at least prove 
that {\it some\/} power of $I$ coincides
with a power of $\mm$, and that in case the resolution
of $I$ is linear for at least $\lceil (n-1)/2\rceil$
steps, then $I^2$ is a power of $\mm$. We can also
give some weak numerical evidence related to the number 
of generators of $I$. This section is devoted to 
these and related more general results.

\theorem{some power} If $I$ is an $\mm$-primary linearly
presented ideal in $n$ variables, generated in degree $d$
$($or, when the ground field is algebraically closed,
if $\mm^d\subset I+(z_3,\dots, z_n)$ for all
sequences of $n-2$ independent linear forms $z_3,\dots, z_n$$)$,
then $I^t=\mm^{dt}$ for all $t\gg 0$.

We will use the following criterion:

\proposition{smoothness} Let $I\subset S$ be an
 ideal generated
by a vector space $V\subset S_d$, for some $d$. If
$I^s=\mm^{ds}$ for some $s$, then 
$I^t=\mm^{dt}$ for all $t\geq s$. 
This condition is satisfied for some $s$ if and only if
the linear
series $|V|$ maps $\PP^{n-1}$ isomorphically to 
its image in $\PP(V)$.

\noindent{\sl Proof of \ref{smoothness}}.
To prove the first assertion it suffices,
by induction, to treat the case $t=s+1$.
Suppose that $I^s=\mm^{ds}$.
Since $I\subset \mm^d$
we get $I\mm^{d(s-1)}=\mm^{ds}$. Thus 
$I^{s+1}= II^s=I\mm^{ds}=I\mm^{d(s-1)}\mm^d= \mm^{ds}\mm^d=\mm^{d(s+1)}$,
as required.

To prove the last assertion, note that the image
of $\PP^{n-1}$ under the map $\phi$
defined
by the linear series $|V|$ is by 
definition the variety with homogeneous
coordinate ring 
$$
\bigoplus_t(V)^{t}\subset \bigoplus_t S_{dt}.
$$
To say that $\phi$ is an isomorphism 
onto its image means that these
two rings are equal in high degree; that is, $(V)^t=S_{d}$,
so $I^t=\mm^{dt}$
for large $t$.\Box 

\smallskip

\noindent{\sl Proof of \ref{some power}}.
We can harmlessly extend the ground field and assume
that it is algebraically closed. 

By \ref{smoothness} it suffices to show that the map
$\phi$ defined by the linear series $|V|$ is an isomorphism.
For this it is even enough to show that the restriction
of $\phi$ to any line is an isomorphism: 
There is a line through any two points of 
$\PP^{n-1}$ and a line containing any tangent vector
to a point of $\PP^n$, so if $\phi$ restricts to an
isomorphism on each line then $\phi$ is one-to-one and
unramified, whence an isomorphism.

A line $\ell\subset \PP^{n-1}$ is defined by an
ideal generated by the vanishing of $n-2$ linear forms, say 
$z_3,\dots,z_n$. The restriction $\phi\mid_\ell$ of $\phi$ to
$\ell$ is defined by the degree $d$ component of the 
ideal $\mm^d\subset I+(z_3,\dots, z_n)/(z_3,\dots,z_n)$.
By \ref{linear complements},
$\mm^d\subset I+(z_3,\dots, z_n)/(z_3,\dots,z_n)=$,
so $\phi\mid_\ell$ is defined by the complete linear series
of degree $d$, which is  an 
isomorphism as required.
\Box

To give the results about \ref{Us} in their
natural generality, we turn to results on the 
regularity of the product of two ideals.

The following result was proved (in a superficially more
special case) by Jessica Sidman [2002].

\theorem{reg of products} Suppose $I$ and $J$ are
homogeneous ideals of $S$ and set $\delta = \dim \Tor_1(S/I,S/J).$
If $j\geq \delta-1$, then
$$
\reg H^{j}_\mm (IJ)\leq \reg I+\reg J.
$$
Thus if $\delta\leq 1$ then
$\reg IJ \leq \reg I+\reg J$, 
and if $\delta\leq 2$ then
$\reg(IJ)^{\rm sat} \leq \reg I+\reg J$.

Since $\Tor_1(A,B)=(I\cap J)/IJ$,
the condition $\dim \Tor_1(A,B)\leq 1$ 
of \ref{reg of Tor} may then be interpreted as
saying that the codimension
of $(IJ)_d$ in $(I\cap J)_d$ is bounded independently
of $d$. Thinking of $I, J$ as determining projective
schemes $X, Y\subset \PP^{n-1}$, we may also state the condition as saying
that $X$ and $Y$ are {\it homologically transverse\/} except
at a finite set of points of $\PP^{n-1}$.

\proof Extending the ground field if necessary,
we may assume it is infinite. 
A general linear form
is then annihilated only by an ideal of finite
length modulo $I,J$ or $IJ$. If $\delta\geq 2$
then factoring out such a general
form, the left hand side of the displayed inequality
can only increase and 
the right hand side can only decrease. Thus it suffices
to treat the case $\delta\leq 1$.

Consider the exact sequence
$$
0\to IJ\to I\to I/IJ\to 0.
$$
Note that $I/IJ=\Tor_0(I,S/J)$.
By \ref{reg of Tor}, $\reg \Tor_0(I, S/J)\leq \reg I+\reg S/J$,
and therefore
\goodbreak
$$\eqalign{
\reg IJ&\leq \max\{\reg I, \reg I/IJ +1\}\cr
&=\max\{\reg I, \reg I+\reg S/J+1\}\cr
&=\reg I+\reg J.
}
$$
\Box

\theorem{reg of products2A}
Suppose that $I$ and $J$ are homogeneous ideals of $S$ with
$\dim \Tor_1(S/I, S/J)\leq 1$. If $p,q$ are integers such that
$p\leq \codim I,\ q\leq \codim J,$ and $p+q=n+1$, then
$$
\reg IJ \leq \max\{\reg I,\ \reg J,\ t_p(S/I)+t_q(S/J)-n+1\}.
$$ 

\proof
{}From the short exact sequences
$$\eqalign{
0\to (I\cap J)/IJ \to &S/IJ \to S/(I\cap J)\to 0\cr
0\to S/(I\cap J)\to S/I&\oplus S/J \to S/(I+J)\to 0
}
$$
we see that 
$$\eqalign{
\reg S/IJ &\leq \max\, \{ \reg S/(I\cap J),\ \reg (I\cap J)/ IJ\}\cr
&\leq
\max\, \{\reg S/I,\ \reg S/J,\ 1+\reg S/(I+J),\ \reg (I\cap J)/ IJ\}.
}$$
Notice that $S/(I+J) = \Tor_0(S/I,S/J)$ and $(I\cap J)/IJ = \Tor_1(S/I,S/J)$.
To bound the regularity of these modules we apply \ref{newregtor}
with $0\leq k\leq 1$. 

From the hypothesis we see that $1+ \dim S/J \leq p\leq \codim I$. Hence by 
\ref{newregtor},
$$
\reg \Tor_0(S/I,S/J)\leq t_p(S/I)+t_{q-1}(S/J)-n
$$
 and 
$$
\reg \Tor_1(S/I,S/J)\leq t_p(S/I)+t_{q}(S/J)-n.
$$
Using the inequalities above, we obtain
$$
\eqalign{
\reg S/IJ \leq \max\, 
\{
\reg S/I,\ 
\reg S/J,\ 
1+&t_p(S/I)+t_{q-1}(S/J)-n,\cr 
  &t_p(S/I)+t_{q}(S/J)-n
\}.
}
$$
Because $q\leq \codim S/J$ we have
$$
 t_{q-1}(S/J)\leq t_q(S/J)-1.
$$
Thus 
$$
\reg S/IJ \leq 
\max\, \{
\reg S/I,\ 
\reg S/J,\ 
t_p(S/I)+t_{q}(S/J)-n\},
$$
as required.\Box

\corollary{reg of products2}
Suppose that $I$ and $J$ are homogeneous ideals of $S$.
If either
$\dim S/J = 0$ and $I$ 
is generated in degrees at most $d$,
or $\dim S/J= 1$ and $I$ 
is related in degrees at most $d+1$, then
$$
\reg IJ \leq \max\, \{\reg I, d+ \reg J\}.
$$

\proof 
We may assume $I\neq 0$, and dividing $I$ by its 
greatest common divisor we may then suppose that 
$\codim I \geq 2$.
We apply \ref{reg of products2A} with $p=1$ in 
the first case, and $p=2$
in the second case.
\Box

\corollary{half-way linear2} Suppose $I$ and $J$ are homogeneous ideals
in $S$ of dimension $\leq 1$, generated in degree $d$.
If the resolutions of $I$ and $J$ are linear for $\lceil (n-1)/2\rceil$ steps
$($for instance if $I$ and $J$ have linear presentation 
and $n\leq 3$$)$, then
$IJ$ has linear resolution.
In particular, if $I$ and $J$ are $\mm$-primary then $IJ = \mm^{2d}$.

\proof  
Applying \ref{newregtor} with $k=0$ we get
$$
\reg S/I=\reg\Tor_0(S/I,S/I)\leq 2d-2,
$$
 and similarly
for $\reg S/J$.
From \ref{reg of products2A} with $p=q=\lceil (n+1)/2\rceil$
we see that $\reg IJ\leq 2d$. Since $IJ$ is generated 
in degree $2d$, it follows that $IJ$ has linear resolution.\Box

Taking $I=J$ we get
the special case $s=\lceil (n-1)/2\rceil$
of \ref{Us}.

\corollary{half-way linear} Suppose $I\subset S$ is a homogeneous ideal
of dimension $\leq 1$, generated in degree $d$.
If the resolution of $I$ is linear for $\lceil (n-1)/2\rceil$ steps $($for
instance if $I$ has linear presentation and $n \leq 3$$)$, then
$I^t$ has linear resolution for all $t\geq 2$.
In particular, if $I$ is $\mm$-primary then $I^2 = \mm^{2d}$.

\corollary{reg of powers} Let $I\subset S$ be a homogeneous ideal
of dimension $\leq 1$.
If $I$ is generated in degree $d$ and has linear 
presentation, and if some power of $I$ has a linear free resolution,
then all higher powers of $I$ have linear free resolution.

\proof Suppose that $I^t$ has a linear
resolution. In \ref{reg of products2}
we replace $I$ by $I^t$, and $J$ by $I$.
It follows that
$I^{t+1}$ has regularity $dt$.
As it is generated in degree $dt$ 
it must have linear resolution.\Box

No such result holds for 2-dimensional ideals in 4 variables $a,b,c,d$:
Aldo Conca [2003] has shown
the ideal $I=(ab^r, ac^r, b^{r-1}cd)+bc(b,c)^{r-1}$, with $r>1$ 
has the property that $I^t$ has linear resolution for $t<r$,
while $I^r$ does not even have linear presentation.
See also Sturmfels [2000].

If $I\subset S$ is an $\mm$-primary ideal generated in 
degrees $\leq d$ then $\reg I^t\leq \reg I +(t-1)d$. 
(Reason: Write $e=\reg I$. Since $\mm^e\subset I$, we have 
$\mm^e\subset\mm^{e-d}I$
and thus 
$
\mm^{e+(t-1)d}\subset \mm^{e+(t-2)d}I.
$
Induction on $t$ completes the argument.)
But we can prove a little more. The following result
is also a generalization of \ref{half-way linear}.

\corollary{reg of powers2} Let $I\subset S$ be a homogeneous ideal
and let $t\geq 2$ be an integer.
If
$\dim S/I = 0$ and $I$ 
is generated in degrees at most $d$
or $\dim S/I= 1$ and $I$ 
is related in degrees at most $d+1$, then
$\reg I^t\leq \reg I + (t-1)d$. More generally, for
$1+\dim S/I \leq p\leq \codim I$,
$$
\reg I^t \leq
 t_{p-1}(I)+t_{n-p}(I)-n+(t-2)d+1.
$$

\proof
We use induction on $t\geq 2$.

\ref{newregtor} shows that
$$\eqalign{
\reg S/I &=
\reg \Tor_0(S/I, S/I)\cr
&\leq t_p(S/I) +t_{n-p}(S/I)-n\cr
&<t_p(S/I) +t_{n+1-p}(S/I)-n,
}
$$
where the last inequality holds because $n+1-p\leq \codim I$.
Similarly,
$$\eqalign{
\reg I^{t-1}/I^t &= \reg \Tor_1(S/I, S/I^{t-1}) \cr
&\leq t_p(S/I) +t_{n+1-p}(S/I^{t-1})-n.
}
$$
Hence the exact sequence 
$$
0\to I^{t-1}/I^t\to S/I^t\to S/I^{t-1}\to 0
$$
shows that
$$\eqalign{
\reg S/I^t &\leq  t_p(S/I)+t_{n+1-p}(S/I^{t-1})-n \cr
&\leq t_p(S/I) +\reg S/I^{t-1} +1-p.
}
$$
The base case $t=2$ of the present corollary now follows from the 
first inequality. The induction step uses the second equality 
with $p=1$ or $p=2$, depending on whether
$\dim S/I = 0$ and $I$ 
is generated in degrees $\leq d$
or $\dim S/I= 1$ and $I$ 
is related in degrees $\leq d+1$.
\Box

There has been considerable recent progress on the general subject of
regularity bounds for powers of an ideal; see Herzog, Hoa and Trung
[2002] and the references cited there.


By
comparing the number of generators of
$\mm^{(n-1)d}$ with the number of generators of the $(n-1)$ symmetric
power of $I$, we see that
\ref{Eisenbud-Ulrich} implies that the minimal number of generators $\mu(I)$
is at least $(n-1)d+1$.  This is exactly the number of generators
of 
$
(x_1,x_2)^{d-1}\mm
$
(\ref{gin conjecture 1} would provide a more precise
version.)

The following
Proposition, when combined with
\ref{linear complements}, 
provides further numerical evidence in the case $d=2$.

\proposition{} Let $I\subset S$
be an $\mm$-primary 
ideal generated by $\mu$ forms of degree $d$. If $\mm^d\subset I+L$
for every ideal $L$ generated by $n-q-1$ independent linear forms,
then 
$$
\mu\geq (q+1)(n-q-1)+{q+d\choose d}.
$$
For example, if $q=1, n=3$ then $\mu\geq d+3$, while if $q=1, d=2$
then $\mu\geq 2n-1$.

\proof
Let $W=S_1$ be the vector space of linear forms in $S$,
and let $V=I_d\subset S_d$.
Consider natural composite map of vector bundles on the Grassmannian
$G$ of
$n-q-1$ dimensional subspaces $\Lambda$
$$
V\to \Sym_d(W) \to \Sym_d(W/\Lambda).
$$
The hypothesis implies that this map is locally everywhere
surjective. 
Because $\Sym_d(W/\Lambda)$ is ample (see Hartshorne [1970, Ch. 3])
the
theorem of Fulton and Lazarsfeld [1981] (see Arbarello,
Cornalba, Griffiths, and Harris [1985], Proposition VII.1.3)
requires that $\dim G\leq \dim V - \rank
\Sym_d(W/\Lambda) +1$, which is the desired inequality.
\Box

We finish this section with a remark about Rees algebras and 
reduction numbers. Recall that if $J\subset I$ are ideals of
$S$, then the {\it reduction number } ${\rm r}_J(I)$ of
$I$ with respect to $J$ is the smallest integer 
$0\leq r \leq \infty$ with $I^{r+1}=JI^r$.

\corollary{Rees}
Let $I \subset S$ be a homogeneous $\mm$-primary ideal 
generated in degree $d$ and assume that $I \neq \mm^d$.
\item{$(a)$} If $I$ has linear presentation, then
$\depth \R(I) =1$.
\item{$(b)$} If the resolution of $I$ is linear for
$\lceil (n-1)/2\rceil$ steps, then 
${\rm r}_J(I)=\max\{2, n-1-\lfloor (n-1)/d\rfloor \}$ 
for every $\mm$-primary ideal $J \subset I$ generated by 
$n$ forms of degree $d$. 
 
\proof
$(a)$ Consider the exact sequence of finitely generated
$\R(I)$-modules
$$
0 \rTo \R(I) \rTo \R(\mm^d) \rTo C \rTo 0.
$$
The module $C \neq 0$ has finite length by \ref{some power},
showing that $\depth \R(I) =1$.

$(b)$ Since $\R(I)$ is not Cohen-Macaulay and $n \geq 2$, one has
${\rm r}_J(I) \geq 2$ according to 
Valabrega-Valla [1978] and Goto-Shimoda [1982]. 
On the hand $I^t=\mm^{dt}$
for every $t\geq 2$ by \ref{half-way linear}. Therefore
$$
{\rm r}_J(I)=\max\{2,{\rm r}_J(\mm^d)\}.
$$

It remains to see that ${\rm r}_J(\mm^d)=e:= n-1- \lfloor (n-1)/d\rfloor$.
As $\reg S/J =n(d-1)$ it follows that $\mm^{de} \not \subset J$, whereas
$\mm^{d(e+1)}\subset J$ and hence $\mm^{d(e+1)}= J \mm^{de}$. Thus indeed
${\rm r}_J(\mm^d)=e$.
\Box

\section{monomial ideals} Monomial ideals
 
In this section we will prove the second statement
of \ref{Us} for monomial ideals,
and give a necessary and sufficient condition for a monomial
ideal to satisfy the asymptotic version.

\theorem{monomial linear} Let $I$ be a linearly presented, $\mm$-primary
monomial ideal in $S=K[x_1,\dots,x_n]$, generated in degree $d$.
If the minimal resolution of $I$ is linear for $s$ 
steps then $I^{t}=\mm^{td}$ for all $t\geq (n-1)/s$.

\ref{monomial linear} follows at once from the next two results:

\proposition{contains}	 If $I$ is an $\mm$-primary monomial ideal
that is generated in degree $d$ and has
linear resolution for $q$ steps, then $I$ contains the ideal
$$
J(d,q)=\sum_{i_1<\cdots<i_{q+1}}(x_{i_1},\dots,x_{i_{q+1}})^d.
$$

\proof Since $I$ is its own initial ideal, in 
any monomial order, the  statement follows from
\ref{initials}.\Box

\proposition{monomial powers} For all $i\geq 1$, 
$J(id,iq)\subseteq J(d,q)^{i}$. In particular,
if $e \geq {n-1\over q}$,
then $J(d,q)^e = \mm^{de}$.

\proof
The second statement follows from the first
because $J(d,q) = \mm^d$ for $q\geq n-1$.

By induction on $i$,
it suffices to show that
$$
J(i\,d,i\,q)\subseteq J(d,q)\;\cdot\; J((i-1)d,(i-1)q).
$$
To this end, let $m=\prod x_j^{a_j}\in J(id,iq)$ be a monomial
of degree $id$. 
By the definition of $J(id,iq)$, at most $iq+1$ of the $a_j$
are nonzero. To simplify the notation we assume that
$a_j=0$ for $j>iq+1.$

Not  every sum  of $q$ of the $a_1,\dots, a_{iq+1}$ 
can be strictly bigger than
$d$; otherwise $id = \sum_j a_j \geq (d+1)i$, a contradiction. Choose
$q$ of the $a_j$ whose sum $\sigma$ 
is maximal with respect to being at most $d$.
By relabeling we may assume these are $a_1,...,a_{q}$. 

Suppose first that there is no
index $k>q$
such that $\sigma + a_k\geq d$.
It follows from the maximality of $\sigma$,
that $a_k\leq a_j$ whenever $j\leq q<k$.
From this we see that the sum of  any  $q+1$ of the $a_j$ is
at most $d-1$. But then 
$$
id = \sum_1^{iq+1}a_j 
\leq (d-1)\lceil{iq+1\over q+1}\rceil \leq (d-1)i,
$$
a contradiction.

Thus there exists an index $k>q$
such that $\sigma + a_k\geq d$. It follows that 
$u:= x_1^{a_1}\cdots x_{q}^{a_{q}}x_k^{d-\sigma}\in J(d,q)$, while
$v := m/u\in J((i-1)d,(i-1)q)$, as required. 
\Box

Here is a criterion for the asymptotic version to hold.

\proposition{monomial criterion} An $\mm$-primary
 monomial ideal $I\subset S$
generated in degree $d$
has a power equal to a power of  $\mm$ if and only if
$J:=\mm (x_1^{d-1},\dots,x_n^{d-1})\subset I$. Further,
$J^e=\mm^{de}$ if and only if $e\geq (d-2)(n-1)$.

The second statement is
proven in the course of the proof of Herzog and Hibi [2003, Theorem 1.1] 
(the original
formulation is for any $\mm$-primary ideal $\mm J'$ 
with $J'$ generated
in degree $d-1$). We include a proof for the reader's
convenience.

\proof Let $V$ be the vector space generated by
the degree $d$ monomials in $I$.
By \ref{smoothness}, $I^e=\mm^{de}$ for some $e$
if and only if the map $\phi$ defined by
 $|V|$ defines an isomorphism. Since 
everything is torus invariant, this is true if and only
if it is true at the fixed points of the torus action. At 
such a fixed point, all but one variable vanishes, say
$x_1=\cdots=x_{n-1}=0$, and $I$ must generate the local ring
of $\PP^{n-1}$ at this point. Thus $I$ must contain $x_ix_n^{d-1}$
for each $i<n$. Since $I$ is $\mm$-primary, it also contains
$x_n^d$, proving the first statement.

Next consider $J=\mm J'$, where $J'=(x_1^{d-1},\dots,x_n^{d-1})$.
The $e$-th power of $J'$ has resolution obtained from
that of the $e$-th power of $\mm$ by substituting $x_i^{d-1}$ 
for $x_i$. Thus the regularity of $J'^e$ is precisely
$(d-1)(e+n-1)-n$, so $J'^e$ contains
$\mm^{(d-1)(e+n-1)-n+1}$ but no lower power.
 Since the generators
of $J'^e$ have degree $(d-1)e$,
we see that $J^e=\mm^eJ'^e=\mm^{de}$ if and only
if $e \geq (d-1)(e+n-1)-n+1-(d-1)e$, that is,
$e\geq (d-2)(n-1)$.
\Box

\section{symmetric algebra section} Torsion in symmetric and
exterior products

In general it seems a difficult problem to understand the
relations defining the {\it Rees algebra\/}
$\R(I):= S\oplus I\oplus I^2\oplus \cdots$
of an ideal $I\subset S$. As a start, we may write $\R(I)$ as a
homomorphic image $\Sym(I)/\A$ 
of the symmetric algebra $\Sym(I)$.
The relations defining $\Sym(I)$ are easily derived from
the relations defining $I$: if $G_1\to G_0\to I\to 0$ is a 
free presentation, then $\Sym(I)=\Sym(G_0)/G_1\Sym(G_0)$.
That is, the defining ideal of $\Sym(I)$ in the polynomial
ring $\Sym(G_0)$ is generated by the image of $G_1$, 
regarded as a space of forms that are linear in the
variables corresponding to generators of $G_0$. 

Thus
the problem is to understand $\A$. Let $\A_t$ be the
component of $\A$ in $\Sym_t(I)$, so that $\A=\oplus_{t\geq 2} \A_t$.
It is easy to see that $\A_t$ is the torsion submodule of $\Sym_t(I)$.
In this section we will study the regularity of 
$\A_t$ in the case where $I$ is a homogeneous
$\mm$-primary ideal. 

An ideal $I$ is said to be of {\it linear type\/} if $\A=0$.
Following Herzog, Hibi and Vladoiu [2003] we say more generally
that $I$ is of {\it fiber type\/}
if $\mm \Sym(I) \cap \A = \mm\A$ or, equivalently,
if a generating set of relations of the {\it fiber ring\/}
$\R(I)/\mm\R(I)$ lifts to a  generating set for $\A$. If
$I$ is generated by forms of degree $d$, then all the 
generators of $\A_t$ have degrees $\geq td$. The simplest
situation occurs when the regularity of $\A_t$ is $td$.

\theorem{partial annihilation} Let $I\subset S$ be a homogeneous
$\mm$-primary ideal.
\item{$($a$)$} If $I$ is generated in degrees at most $d$ 
and related in degrees at most $e+1$, then 
$\reg \A_t\leq \reg I + (t-2)d+e$ for every $t$.
\item{$($b$)$} Suppose that $I$ is generated
in degree $d$ and has linear presentation. Let $s$ 
be an integer such that 
$I^s=\mm^{sd}$. We have
$\reg \A_{s+u}\leq \max\{ \reg \A_s, sd\} + ud$ for every $u\geq 0$.
\item{$($c$)$} If the resolution of $I$ is linear for
$\lceil n/2\rceil$ linear steps, then 
$\A_t$ is concentrated in degree $dt$ for every $t$;
in particular, $I$ is of fiber type and $\A$
is annihilated by $\mm$.

In the course of their study of implicitization of surface,
Bus\'e and Jouanolou [2003, Prop 5.5] proved  a different bound for the 
torsion in the symmetric
algebra $\Sym(I)$ for ideals $I$
of dimension $\leq 1$.  
This was later sharpened by Bus\'e and Chardin
[2004]. (Although the result was originally
stated only for ideals with $n+1$ generators,
this restriction is irrelevant. A forthcoming paper of
Chardin will contain further generalizations.)

Our proof of \ref{partial annihilation}
is based on a more general lemma:

\lemma{reg of A} If $I\subset S$ is a homogeneous
$\mm$-primary ideal generated in degrees
at most $d$ then 
$$
\reg \A_{t+1}\leq \max\, \{d+\reg \A_t,\ 
\reg \Tor_2(S/I, S/I^t)\}.
$$

\noindent{\sl Proof of \ref{reg of A}}.
Let $G_1\to G_0\to I$ be a minimal free
presentation, so that $G_0$ is generated
in degrees $\leq d$. There is an commutative diagram
with exact rows and columns of the form
$$\diagram[small, tight]
&&&&&&0&& \cr
&&&&&&\dTo&& \cr
&&&&&&\Tor_2(S/I,S/I^t)&& \cr
&&&&&&\dTo&& \cr
&&G_0\otimes\A_t&\rTo &I\otimes \Sym_t(I)&\rTo&I\otimes I^t&\rTo&0\cr
&&\dTo&&\dTo&&\dTo&& \cr
0&\rTo&\A_{t+1}&\rTo&\Sym_{t+1}(I)&\rTo& I^{t+1}&\rTo&0\cr
&&&&\dTo&&\dTo&& \cr
&&&&0&&0&& \cr
\enddiagram
$$
where the left-hand map is given by the $\Sym(G_0)$-module
structure on $\Sym(I)$. The Snake Lemma shows that
$\A_{t+1}$ is an extension of a quotient of $G_0\otimes \A_t$
by a quotient of $\Tor_2(S/I, S/I^t)$. Since both these
modules have finite length, the regularity of such an 
extension is bounded by the maximum of the two regularities
as required.\Box

\smallskip
\noindent{\sl Proof of \ref{partial annihilation}}.
(a) We do induction on $t$. If $t\leq 1$ then 
$\A_t=0$ so the assertion is trivial.
For $t\geq 2$ we apply \ref{reg of A}, and
it suffices to prove 
$\reg \Tor_2(S/I,S/I^{t-1})\leq \reg I+ (t-2)d+e$.
From \ref{tor-bound1} with $p=2$ we obtain
$$
\reg \Tor_2(S/I, S/I^{t-1})
\leq t_2(S/I)+t_n(S/I^{t-1})-n
\leq e+1+\reg S/I^{t-1}.
$$
Hence 
$\reg \Tor_2(S/I,S/I^{t-1})\leq \reg I^{t-1}+e\leq \reg I+ (t-2)d+e$
by \ref{reg of products}, as required.

(b) The same argument works, but this time we start the
induction from $t=s$, and use the fact that 
$\reg I^{s+u-1}=(s+u-1)d$ by \ref{smoothness}.

(c) By \ref{half-way linear} we know that $I^2=\mm^{2d}$, and from
\ref{reg of A} we have
$\reg \A_2\leq \reg\Tor_2(S/I,S/I)$. By
\ref{tor-bound1} with $p=\lceil n/2\rceil +1$
we obtain
$\reg\Tor_2(S/I,S/I)\leq 2d$. Thus
we can apply part (b) with $s=2$ to obtain the desired
result.\Box

\example{non annihilation}
The conclusion of \ref{partial annihilation} $(c)$
does not hold for linearly presented monomial ideals
in 3 variables. For example, let $I$ be the ideal in 
$K[x,y,z]$ generated
by all the monomials of degree 5 except $x^3yz, xy^3z, xyz^3$.
The ideal $I$ is linearly presented, but (since it is not $\mm^3$)
it does not have $\lceil n/2\rceil=2$ linear steps in its resolution.
In this case the module $\A_2$
is generated in degree $2d=10$, but has regularity 11
instead of 10. (In this case $\A_t$ does have regularity
$5t$ for all $t\geq 3$.)

The conclusion of \ref{partial annihilation} does hold
for linearly presented primary ideals in 3 variables in the
Gorenstein case:

\corollary{Gorenstein} Let $I\subset K[x_1,x_2,x_3]$ be a homogeneous
$\mm$-primary Gorenstein ideal. If $I$ is generated in degree $d$ and 
has linear presentation, then
$\A_t$ is concentrated in degree $dt$ for every $t$;
in particular, $I$ is of fiber type and $\A$
is annihilated by $\mm$.

\proof
We know that $I^2=m^{2d}$ by \ref{half-way linear} and $\A_2=0$
by Huneke [1984]. Hence the assertion follows from part (b) of
\ref{partial annihilation}.
\Box

The application of \ref{newregtor} to $\Tor_2$ also yields a
result on the regularity of exterior powers:

\corollary{exterior torsion} If $\dim S/I\leq 1$ then
$$
\reg\H^0_\mm(\wedge^2 I)\leq 
\reg \H^0_\mm(I\otimes I)
\leq  t_p(S/I)+ t_q(S/I)-n.
$$
for any $p,q\leq \codim I$ such that $p+q=n+2$. In particular,
if $I$ is an $\mm$-primary ideal generated in degree $d$ with
linear free resolution for $\lceil n/2\rceil$ steps, then
$\wedge^tI$ is a vector space concentrated in degree $dt$
for every $t\geq 2$.

\proof For the first statement we simply observe that
the torsion submodule of $I\otimes I$ is $\Tor_2(S/I,S/I)$,
and apply \ref{newregtor}. 

To obtain the second statement for $t=2$
we apply the first statement with
$p=\lceil (n+2)/2\rceil$ and 
$q=\lfloor (n+2)/2\rfloor$. Since $\wedge^2I$ is
always annihilated by $I$, we have $\H^0_\mm(\wedge^2I) = \wedge^2I$
in this case. For general $t\geq 2$ we use induction,
noting that there is always a surjection 
$\wedge^2I\otimes \wedge^{t-2}I\to \wedge^tI$.
\Box

\section{elimination section} Application: Instant elimination

Let $I$ be an ideal of $S$, generated by a vector space
$V$ of forms of degree $d$. We may think of $V$ as a
linear series on $\PP^{n-1}$ and ask for the equations
of the image scheme; we may also restrict $V$ to a subscheme
$X\subset \PP^{n-1}$ to try to compute the image of $X$.
These computations involve the elimination of 
variables: If $V = \langle f_1,\dots,f_m\rangle$ then
we are looking for the relations on the elements $f_it$ in 
$S_X[It]\subset S_X[t]$. Geometrically, the ideal $I$ defines 
the base locus of a blowup, and we are looking for the defining
relations on the {\it fiber\/}
$\R_{S_X}(I)/\mm\R_{S_X}(I)$.

In some interesting classical cases, there is a much easier way to do 
elimination. For example, if $V$ is the linear series of 
$d$-ics through a set $B$ of ${d+1\choose 2}$ 
general points in the projective
plane
then the ideal $I$ generated by $V$ is linearly
presented: indeed, by the Hilbert-Burch theorem,
 the free resolution of $S/I$ has the form
$$
0\rTo S(-d-1)^{d}\rTo^{\phi} S(-d)^{d+1}\rTo S.
$$
The $d\times (d+1)$ matrix $\phi$ of linear forms in 3 variables 
may be thought of as a $d\times d\times 3$ tensor
over $K$. This tensor may also be identified with
a matrix
$\psi$ of size $3\times d$ in $d+1$ variables, called
the {\it adjoint\/} (or {\it Jacobian dual\/}) matrix. The image of
$\PP^2$ under the rational map defined
by $V$ is isomorphic to $\PP^2$ blown up
at $B$. The defining ideal of this variety is generated by
the $3\times 3$ minors of $\psi$ (Room [1938]) see also
Geramita and Gimigliano [1991], and  Geramita, 
Gimigliano and Pittleloud [1995], which
does the case of determinantal sets of points in $\PP^r$.

The idea of doing elimination in this way was generalized
and put to practical use
by Schreyer and his coworkers 
(Decker-Ein-Schreyer [1993] 
Ranestad-Popescu [1996], Popescu [1998]) in their study of 
surfaces of low degree in $\PP^4$, in cases where the usual
elimination methods were too demanding computationally. 
It is easy to see that the method works whenever $I$ is
of linear type (as an ideal of $S_X$, 
in the sense that the powers of $I$ are
equal to the symmetric powers). But the examples above
are not of linear type. 

Here is a general criterion for when the instant elimination
process works.
We regard $\Sym I$ and $\R(I)$ as bigraded algebras
with degrees with an element of degree $a$ in $\Sym_b(I)$
being given degree $(a,b)$.

\proposition{condition for instant elimination} Let $X\subset \PP^{n-1}$
be a scheme, and let $V$ be a linear series of forms of degree $d$ on $X$.
Suppose that the ideal $I$ generated by $V$ has
linear presentation, with matrix $\phi$, and that $\psi$ is
the adjoint matrix. If the torsion in the symmetric algebra
of $I$ occurs only in degrees $(a,b)$ such that $a=db$,
then the annihilator of $\coker\psi$ is the ideal of 
forms in $\PP(V)$ that vanish on the image of $X$ under the
rational map associated to $|V|$.

\proof Write $V=\langle f_1,\dots,f_m\rangle$. 
We consider the epimorphism of bigraded algebras
$$
K[X_1,\dots,X_n, T_1,\dots,T_m]\to \Sym(I);\quad 
X_i\mapsto x_i, \quad T_i\mapsto f_i\in \Sym_1(I)
$$
where $X_i$ is an indeterminate of degree $(1,0)$ and
$T_i$ is an indeterminate of degree $(d,1)$. 
There are $K[T_1,\dots,T_m]$-module isomorphisms
$$
\coker \psi 
\cong \oplus_b(\Sym(I))_{(bd+1,b)}
\cong \oplus_b(\R(I))_{(bd+1,b)},
$$
where the last isomorphism follows from our assumption
about the torsion of $\Sym(I)$. On the other hand,
since $\R(I)$ is a domain, 
$\sum_d(\R(I))_{(bd+1,b)}$ 
and
$\sum_d(\R(I))_{(bd,b)}=K[f_1t,\dots,f_mt]$
have the same annihilator.
\Box

\corollary{instant elimination} Suppose that $V$ is a
base point free linear
series of forms of degree $d$ on $\PP^{n-1}$. 
Suppose that the free resolution of the ideal $I$ generated by $V$
is linear for at 
least $\lceil n/2\rceil$ steps.
Let $\phi$ be the presentation matrix of $I$. If $\psi$ is the
adjoint matrix of $\phi$ then
the annihilator of $\coker\psi$ is the ideal of 
forms in $\PP(V)$ that vanish on the image of $\PP^{n-1}$ under the
rational map associated to $|V|$.

\proof Apply \ref{partial annihilation} $(c)$ and 
\ref{condition for instant elimination}.\Box

\section{almost lin} Ideals with almost linear resolution

We can get a bound for the number of generators of an
ideal with ``almost linear'' resolution as follows.
Let $n=r+1$ so that $S=K[x_0,\cdots,x_r],$ with $r\geq 2$ to 
avoid the trivial case,
and suppose that the free resolution of $S/I$ has the form
$$
S\lTo S^{m_1}(-d)\lTo\cdots\lTo S^{m_r}(-d-r+1)\lTo
\sum_1^{m_{r+1}}
S(-d-r-b_i)\lTo 0;
$$
that is,  $I$ is generated in degree $d$,
$S/I$ has ``almost linear resolution'', 
and the socle elements of $S/I$ lie in degrees $d+b_i-1$,
with $b_i\geq 0$.
Assume further that $S/I$ has finite length. Our goal
is to find a lower bound for the number of generators
of $I$.

Computing the Hilbert polynomial $0\equiv P_{S/I}(\nu)$ we get
$$
0=
{\nu+r\choose r}
+\sum_{i=1}^r (-1)^i m_i {\nu - d - (i-1)+r\choose r}
+(-1)^{r+1}\sum_1^{m_{r+1}} {\nu-d-b_i\choose r}.
$$
Taking $\nu=d-1$, all but the first and last terms vanish, so
$$
{d+r-1\choose r}=(-1)^r\sum{-b_i-1\choose r}=\sum {b_i+r\choose r}.\eqno{(1)}
$$
Taking $\nu=d$, all but the first two and the last terms vanish, so
$$
m_1={d+r\choose r}-\sum {b_i+r-1\choose r},\eqno{(2)}
$$
or equivalently $\H_{S/I}(d)= \sum {b_i+r-1\choose r}$.

Continuing in this way we could inductively compute all the 
$m_i$ in terms of the $b_i$. But already equations $(1)$
and $(2)$ suffice to give a lower bound for the number
of generators:

\proposition{generator bound} With notation as above,
$$
m_1\geq {d+r-1\choose r-1}+{d+r-2\choose r-1}
$$
with equality if and only if $S/I$ is Gorenstein.

\proof
By equation $(1)$ we have $d-1\geq b_i$ for every $i$, and equality
holds for some $i$ if and only if $m_{r+1}=1$, 
that is, if $S/I$ is Gorenstein (and there is only one $b_i$.)
Thus by equation $(2)$
$$\eqalign{
m_1&={d+r\choose r}-\sum {b_i+r-1\choose r}\cr
&=
{d+r\choose r}-\sum {b_i\over b_i+r}{b_i+r\choose r}\cr
&\geq
{d+r\choose r}-{d-1\over d+r-1}\sum {b_i+r\choose r}}
$$ 
with equality if and only if $S/I$ is Gorenstein.
By equation $(1)$ we may rewrite the last line 
as 
$$
{d+r\choose r}-{d-1\over d+r-1}{d+r-1\choose r}
={d+r-1\choose r-1}+{d+r-2\choose r-1}.\Box
$$

\references

\noindent 
E.~Arbarello, M.~Cornalba, P.~A.~Griffiths and J.~Harris:
{\sl Geometry of algebraic curves\/} Vol. I. 
Grundlehren der Mathematischen Wissenschaften 
267. Springer-Verlag, New York, 1985.

\medskip\noindent
M.~Auslander: Modules over unramified regular local rings. 
Illinois J. Math. 5 (1961) 631--647.

\medskip\noindent
L. Bus\'e and M. Chardin:
Implicitizing rational hypersurfaces using approximation complexes. 
J. Symbolic Computation, to appear. 

\medskip\noindent
L. Bus\'e and J.-P. Jouanolou: 
On the closed image of a rational map and the implicitization problem.
J. Algebra 265 (2003) 312--357.

\medskip\noindent
G. Caviglia:
Bounds on the Castelnuovo-Mumford regularity of tensor products.
Proc. Am. Math. Soc., to appear.

\medskip\noindent
M. Chardin:
Preprint (200?)

\medskip\noindent
A. Conca: Regularity jumps for powers of ideals.\hfill\break
http://www.arxiv.org/math.AC/0310493

\medskip\noindent
S. D. Cutkosky: Irrational asymptotic behaviour of Castelnuovo-Mumford
regularity. J. Reine Angew. Math. (2000) 93--103.

\medskip\noindent 
S. D. Cutkosky, J. Herzog and N. V. Trung:
Asymptotic behaviour of the Castelnuovo-Mumford regularity.
Composition Math. 118 (1999) 243--261.

\medskip\noindent
W. Decker, L. Ein, F.-O.  Schreyer: 
Construction of surfaces in $ P\sb 4$.
J. Algebraic Geom. 2 (1993) 185--237. 

\medskip\noindent
D.~Eisenbud: {\sl Commutative algebra with a view towards
algebraic geometry.\/} Spinger-Verlag, New York, 1995.

\medskip\noindent
D.~Eisenbud: {\sl Geometry of Syzygies.\/} 
Spinger-Verlag, New York, 2004 (in preparation; 
available at www.msri.org/\~{}de/ready.pdf ).

\medskip\noindent
D. Eisenbud, M. Green, K. Hulek and S. Popescu:
Restricting linear syzygies: algebra and geometry.
http://www.arxiv.org/math.AG/0404517 (2004a).

\medskip\noindent
W.~Fulton and R.~Lazarsfeld: 
On the connectedness of degeneracy loci and special
divisors.  Acta Math. 146 (1981) 271--283.


\medskip\noindent
A. V. Geramita and A. Gimigliano:
Generators for the defining ideal of certain rational surfaces.
Duke Math. J. 62 (1991)  61--83.

\medskip\noindent
A. V. Geramita, A. Gimigliano, and Y.  Pitteloud:
Graded Betti numbers of some embedded rational $n$-folds.
Math. Ann. 301 (1995)  363--380.

\medskip\noindent
S. Goto and Y. Shimoda:
On the Rees algebras of Cohen-Macaulay local rings.
Commutative algebra (Fairfax, Va., 1979) pp. 201--231,
Lecture Notes in Pure and Appl. Math., 68,
Dekker, New York, 1982.

\medskip\noindent
T. Harima, J. Migliore, U. Nagel and J. Watanabe:
The Weak and Strong Lefschetz Properties for Artinian K-Algebras.
J.  Algebra, 262 (2003)  99--126.

\medskip\noindent
R.~Hartshorne: 
Ample vector bundles. 
Inst. Hautes ƒtudes Sci. Publ. Math. No. 29 (1966) 63--94.

\medskip\noindent
R.~Hartshorne: 
{\sl Ample subvarieties of algebraic varieties. \/}
Notes written in collaboration with C. Musili. 
Lect.~Notes in Mathematics, Vol. 156, 
Springer-Verlag, Berlin-New York 1970.

\medskip\noindent
J. Herzog and T. Hibi:
Castelnuovo-Mumford regularity of simplicial semigroup rings with isolated singularity.
Proc. Amer. Math. Soc. 131 (2003) 2641--2647.

\medskip\noindent
J. Herzog, T. Hibi and M. Vladoiu: 
Ideals of fiber type and polymatroids. 
http://www.uni-essen.de/\~{}mat306/preprints/StrongPol.ps .
Preprint (2003). 

\medskip\noindent
J. Herzog and N.V. Trung: 
Asymptotic linear bounds for the Castelnuovo-Mumford regularity.
Trans. Amer. Math. Soc. 354 (2002) 1793--1809.

\medskip\noindent
J. Herzog, L.T. Hoa and N.V. Trung: Asymptotic linear bounds for the
Castelnuovo-Mumford regularity.
Trans. Amer. Math. Soc. 354 (2002) 1793--1809.

\medskip\noindent
C. Huneke:
Numerical invariants of liaison classes.
Invent. Math. 75 (1984) 301--325.

\medskip\noindent
V. Kodiyalam:
Asymptotic behaviour of Castelnuovo-Mumford regularity. 
Proc. Amer. Math. Soc. 128 (2000) 407--411.


\medskip\noindent
S. Popescu, K.  Ranestad: 
Surfaces of degree $10$ in the projective fourspace via linear systems and linkage.
J. Algebraic Geom. 5 (1996)  13--76.

\medskip\noindent
S. Popescu: 
Examples of smooth non-general type surfaces in $ P\sp 4$.
Proc. London Math. Soc. (3) 76 (1998)  257--275.

\medskip\noindent
T. G. Room:
The geometry of determinantal loci.
Cambridge [Eng.] University Press, 1938.

\medskip\noindent
J. Sidman:
On the Castelnuovo-Mumford regularity of products of ideal.
sheaves. Adv. in Geometry 2 (2002) 219--229.

\medskip\noindent
B. Sturmfels: Four counterexamples in combinatorial algebraic geometry.
J. Algebra 230 (2000) 282--294.

\medskip\noindent
P. Valabrega and G. Valla:
Form rings and regular sequences.
Nagoya Math. J. 72 (1978) 93--101.

\bigskip
\noindent David Eisenbud\hfill\break
Mathematical Sciences Research Institute\hfill\break
17 Gauss Way\hfill\break
Berkeley CA 94720\hfill\break
de@msri.org\hfill\break
\smallskip
\noindent Craig Huneke\hfill\break
Dept. of Mathematics\hfill\break
University of Kansas\hfill\break
Lawrence, KS 66045\hfill\break
huneke@math.ukans.edu\hfill\break
\smallskip
\noindent Bernd Ulrich\hfill\break
Dept. of Mathematics\hfill\break
Purdue University \hfill\break
W. Lafayette, IN 47907\hfill\break
ulrich@math.purdue.edu\hfill\break
\end